\documentclass[a4,twoside,11pt,leqno,openright]{article}
\usepackage[utf8]{inputenc}	
\usepackage[english]{babel}
\usepackage{latexsym}
\usepackage[T1]{fontenc}
\usepackage{enumerate}
\usepackage{a4wide}
\usepackage{fancyhdr}
\usepackage{amsmath}
\usepackage{amsfonts}
\usepackage{amssymb}
\usepackage[all]{xy}
\usepackage{amscd}
\usepackage{amstext}
\usepackage{theorem}
\newtheorem{theo}{Theorem}[section]
\newtheorem{defi}{Definition}[section]
\newtheorem{lem}{Lemma}[section]
\newtheorem{prop}{Proposition}[section]
\newtheorem{coro}{Corollary}[section]
\newenvironment{Proof of Theorem 1.1}%
{\noindent{\textit{Proof of Theorem 1.1~:}}\small}%
{\normalsize\rmfamily\\[.2cm]} 
\newenvironment{Proof of Theorem 1.2}%
{\noindent{\textit{Proof of Theorem 1.2~:}}\small}%
{\normalsize\rmfamily\\[.2cm]}
\newenvironment{preuve}
{\noindent{\textit{Proof~:}}\small}%
{\normalsize\rmfamily\\[.2cm]}
\newenvironment{rem}
{\noindent{\textbf{Remark~:}}} %
{\normalsize\rmfamily\\[.2cm]}
\begin{document}
%\begin{center}
%	{\large {\bf Lacunarity and cyclic vectors for the Backward shift}}\\
%\vspace{0.3cm}
	%{\large {\bf DANS LES ESPACES DE HARDY  VALEURS VECTORIELLES.}}
%\end{center}
%\vspace{1cm}
\title{Lacunarity and cyclic vectors for the Backward Shift}
\date{}
\author{R\'{e}da Choukrallah}
\maketitle
\begin{abstract}
This article gives a description of invariant subspaces for the backward shift generated by vector valued lacunary series and by a class of lacunary power series in $H^2(\mathbb{D},\, X)$, (where $X$ is an Hilbert space). In particular, we show that these series $f$ in $H^2(\mathbb{D},\, X)$ are cyclic vectors if and only if the queue of Taylor coefficients $\{\widehat{f}(k)$, $k>N\}$ generates the whole space $X$. Analogues of this result are obtained for some functions whose spectrum is a finite union of lacunary sequences and in the polydisc. In the scalar case $H^2$, we give a criterion on the Fourier spectrum of the function to have cyclicity for any power of the backward shift.
\end{abstract}
%\title{Lacunarity and cyclic vectors for the backward shift}
%\maketitle
%\newpage
%\thispagestyle{empty}
%\setcounter{page}{1} 
%\tableofcontents
%\newpage
%\thispagestyle{empty}
\section*{Introduction}
%\addcontentsline{toc}{chapter}{Introduction}
%\markright{Introduction}
$H^2(\mathbb{D},\, X)$ is the Hardy space of $X$-valued functions $f$ that are holomorphic in the disc $\mathbb{D}=\{ \zeta:\, \zeta \in \mathbb{C}\, , \,| \zeta |<1 \}$ and such that,\\ 
$$\| f \|^2:=\displaystyle{\sup_{0\le r<1}} \displaystyle{\int_{ \mathbb{T}}}\|f(r \zeta)\|_{X}^{2} dm(\zeta)< \infty \, ,$$
where $\mathbb{T}=\{ \zeta:\,  | \zeta |=1 \}$ is the unit circle and $m$
is the normalized Lebesgue measure on $\mathbb{T}$. We will write $H^2$ for $H^2(\mathbb{D},\, \mathbb{C})$.\\
The shift operator $S$ and its adjoint $S^{*}$ act in $H^2(\mathbb{D},\, X)$ by the formulas,
$$S_Xf=zf \, , \, \, \, S_X^{*}f=\frac{f-f(0)}{z} \, ,$$
where $z$ is the ``independant variable'' i.e. the identity mapping of the disc $\mathbb{D}$ (or circle $\mathbb{T}$) onto itself ($z( \zeta) \equiv \zeta).$ For shortness and when it is clear in which spaces these operators are acting, we will write $S$ and $S^*$ for $S_X$ and $S_X^*$.\\
It is well known (see $\cite{FoNa}$) that $H^2(\mathbb{D},\, X)$ can be described as the space $\ell^2_a(X)$ of power series $f(z)= \displaystyle{\sum_{n \ge 0}} \hat{f}(n) z^n$ such that $ \hat{f}(n) \in X \,$ and $\, \displaystyle{\sum_{n \ge 0}} \| {\hat{f}(n) \|}_X^2 < \infty$ and if a function $g \in H^2(\mathbb{D},\, X)$ is represented by the sequence of its Fourier coefficients, $g=\{ \widehat{g}(0), \, \widehat{g}(1), \, \widehat{g}(2), \, \ldots \}$, $S$ and $S^*$ are called respectively the forward (right)and the backward (left) shifts,
$$S \{ \widehat{g}(0), \, \widehat{g}(1), \, \widehat{g}(2), \, \ldots \}=\{0,\, \widehat{g}(0), \, \widehat{g}(1), \, \ldots \}.$$
$$S^* \{ \widehat{g}(0), \, \widehat{g}(1), \, \widehat{g}(2), \, \ldots \}=\{ \widehat{g}(1), \, \widehat{g}(2), \, \widehat{g}(3), \, \ldots \}.$$ 
Giving a familly of functions $F\subset H^2(\mathbb{D},\, X)$, we consider the $S ^*$-invariant subspace
$$E_F \stackrel{def}{=}span(S^{*n} f:n \ge  0,\, f \in F).$$ 
A familly of functions $F\subset H^2(\mathbb{D},\, X)$ is said to be cyclic for the backward shift if 
$$E_F=H^2(\mathbb{D},\, X).$$
In the case of a finite dimensional space $X=\mathbb{C}^d$, the space $H^2(\mathbb{D},\, X)$ can also be described in others ways: in terms of coordinate functions, namely, $f=(f_{i})_{1 \le i \le d} \, , \, f_{i} \in H^2$ and $\|f \|^2= \displaystyle{\sum_{1}^{d}} \| f_i \|^2$ and we can represent the shift operateurs as an orthogonal sum of scalar shifts,
$$S_d=S\oplus \ldots \oplus S \, : \, H_d^2 \longrightarrow H_d^2 \, \, \mathrm{et} \, \,  S_d^*=S^*\oplus \ldots \oplus S^* \, : \, H_d^2 \longrightarrow H_d^2,$$
where $H_d^2=H^2 \oplus \ldots \oplus H^2$. Then, the cyclicity of a vector-valued function $f=(f_{i})_{1 \le i \le d}$ in $H^2(\mathbb{D},\, \mathbb{C}^d)$ means the possibility of "simultaneously approximating" any set $(g_{i})_{1 \le i \le d}$ where $g_{i} \in H^2 \, , \,1 \le i \le d$: there is a sequence of complex polynomial $(p_n)_{n \ge 1}$ such that,
$$\displaystyle{ \lim_{n} p_{n}}(S^*)f_{i}=g_{i} \, ,\, \,1 \le i \le d.$$
The "Fourier spectrum" (called also the frequency spectrum) and for shortness, the spectrum of the function $f \in H^2(\mathbb{D},\, X)$ is the set $$\sigma(f)=\sigma_{\mathcal{F}}(f)=\{ k \ge 0:\hat{f}(k) \ne 0 \}.$$
A function $f$ holomorphic in $\mathbb{D}$ is said to be $\sigma$-$spectrale$ if $\sigma(f) \subset \sigma$ and Hadamard lacunary power series are functions of the form;\\
$$f(z)=\displaystyle{\sum_{k=1}^{\infty}} a_k z^{n_k}, \, \, and \, \, such \, \,that \, \, \, \, \frac{n_{k+1}}{n_k} \ge d>1 \, , \, \, \, \, \forall \, \,k,$$
(the constant $d$ is independant on $k$) then  $f$ is $\sigma$-spectral and $\sigma$ is a lacunary set (in the sense of Hadamard).\\
Given a set $\Lambda \subset \mathbb{N}$, we note
$$H_\Lambda^2(\mathbb{D},\, X)=\{f \in H^2(\mathbb{D},\, X):\, \sigma(f) \subset \Lambda \}.$$
$$\ast \, \ast$$
The problem of cyclicity for an operator $T$ in an Hilbert space is connected to the problem of the existence of non-trivial $T$-invariant subspaces. More precisely, there exists a none cyclic vector $x \ne 0$ if and only if there exists a none trivial $T$-invariant subspace and the basic motivations for the study of invariant subspaces  come from interest in the structure of operators and from approximation theory.\\
This article deals with the phenomen of cyclicity for the backward shift. In the first part, we want to give a description of $S^*$-invariant subspaces generated by a class of lacunary series and we will study the open problem of the cyclicity of lacunary series in $H^2(\mathbb{D},\, X)$ where $X$ is a separable Hilbert space. In particular, we obtain an explicit criterion for the cyclicity of lacunary series which sequence formed by its Taylor coefficients is completly relatively compact (c.r.c.). As we will see forwards in details, a sequence $(x_n)_{n \ge 1} \subset X$ is said to be c.r.c. if for any orthogonal projections $P:\, X\rightarrow X$, the normalized sequence $\{ Px_n / \|Px_n\|:\, Px_n \ne 0 \}$ is relatively compact. We will see that this class of sequences coincides with all the sequences when $X$ is of finite dimension and for every separable space $X$, there exists c.r.c. sequences which generated the whole space $X$. In the second part of this paper, we will see how we can connect our results to the scalar case and give a criterion of cyclicity for any power of the backward shift. The third part deals with some particular series whose spectrum is a finite union of lacunary series and in the fourth part we generalize some of our results to the polydisc.\\
Recall that the first results of cyclicity in the scalar case $H^2$ were obtained by R. Douglas, H. Shapiro and A. Shields in 1970 (see $\cite{DSS}$) which is a reference for the study  of cyclicity of the backward shift and source of inspiration for many autors. In fact, their paper contains two approches concerning the cyclicity of lacunary series.  The first approach is based on a pseudocontinuation and some purely arithmetic properties of the spectrum. By this approach, they proved that if $f(z)= \displaystyle{\sum_{k=0}^\infty} a_k z^{2^k} \in H^2$ and if $a_k \ne 0$ for an infinity of $k$ then $f$ is cyclic for the backward shift. This approach can be applied to some not so spare spectrum as for example sets of the form $E=\{n^2+m^2: \, n, \, m \in \mathbb{Z} \}$ (see \cite{Alek0}) or also $E=\{p_n:\, n \in \mathbb{N} \}$ where $\{p_n \}_{n \in \mathbb{N}}$ is the sequence of prime numbers. On the other hand, it is not clear whether this method can be applied to some irregular Hadamard lacunary series.\\
The second approach is to work with Taylor coefficients and consider relations between the spectrum and the approximation ability of the whole space by linear combinations of the truncated queue of Taylors coefficients $S^{*n}f(z)=\displaystyle{\sum_{k \ge 0}} \widehat{f}(k+n)z^k$. This technique uses the following property verify by lacunary series,
$$\displaystyle{\sup_{I \in \mathbb{N}}} \, card \big \{ (m, \, n) \in \Lambda^2: \, m \ne n ,\, m-n=I \big \} < \infty.$$
R. Douglas, H. Shapiro and A. Shields proved that a lacunary serie which is not a polynomial is cyclic in $H^2$ and E. Abakumov in his paper $\cite{Abak}$ proved the result under the weaker condition that the spectrum is a finite union of lacunary sets. We also recall the result of A. B. Aleksandrov $\cite{Alek}$ who proved that series whose spectrum is infinite and included in a $\Lambda(1)$ set (which is a more general set than lacunary sets, finite union of lacunary sets and even Sidon sets) are cyclic.\\ 
\\
What can we say about cyclicity of lacunary series in the more general Hardy space $H^2(\mathbb{D},\, X)$ with values in an Hilbert space $X$. The problem of cyclicity when the dimension of $X$ is finite was raised by N. K. Nikolskii and V. I. Vasyunin (see \cite{NiVa2}) who introduced a classification and give a study of functions in $H^2(\mathbb{D},\, \mathbb{C}^d)$ according to their degree of non-cyclicity to analyse the phenomen of cyclicity in the vector-valued case. Our approach here is based on the use of Taylor coefficients and the properties on the spectrum of lacunary series.\\ 
\section{Lacunary series}
\subsection{Completly Relatively compact sequences}
\begin{defi}
Let $X$ be an Hilbert space.\\
A sequence $(a_k)_{k \ge 0}$ of elements in $X$ is said completly relatively compact (c.r.c.) if for any orthogonal projection $P:\, X \rightarrow X$ and $\eta =\{k \ge 0:\, Pa_k \ne 0 \}$, the sequence $\big ( \frac{Pa_k}{\| Pa_k \|} \big )_{k \in \eta}$ is relatively compact in $X$.
\end{defi}
\begin{rem} In a finite dimensional space, every sequence is c.r.c.
\end{rem}
In the following Lemma we can see that even in the more general case of Banach spaces of infinite dimension, there exist sequences which generate the whole space and verifying a condition nearly the same as c.r.c. obtened by replacing the orthogonal projections by linear bounded operator. 
\begin{lem}\label{crc}
For any separable Banach space $X$, there exists a sequence $(a_k)_{k \ge 0},\, a_k \in X$ such that 
\begin{enumerate}[1)]
\item $span(a_k:\, k \in A)=X$ for any infinite set $A \subset \mathbb{N}$.
\item For any bounded linear operator $T:\, X \rightarrow X$,
$$\big ( \frac{Ta_k}{\| Ta_k \|} \big )_{k \in \eta}\, ,\, \, \eta =\{k \ge 0:\, Ta_k \ne 0 \}$$
is relatively compact in $X$.
\end{enumerate}
\end{lem}
\begin{preuve}
Let $(x_k)_{k \ge 1} \subset X$ be a normalized sequence such that $span(x_k:\, k\ge  1)=X$ and,
$$a_k= \displaystyle{\sum_{j \ge 0}} \lambda_k^j x_j\, ;\, \, \Bigg \{ \begin{array}{c} \lambda_k \in \mathbb{C}^* \\ |\lambda_k |<1,\, \lambda_i \ne \lambda_j,\, i \ne j \\ \displaystyle{\lim_{k \rightarrow + \infty}} \lambda_k=0. \end{array}$$
Then $a_k=f(\lambda_k)\, \, \, \forall \, \, k \ge 1$ where $f(z)= \displaystyle{\sum_{j \ge 1}}z^jx_j$ is holomorphic in $\mathbb{D}$ with values in $X$.\\
To prove that the sequence $(a_k)_{k \in A}$ generates $X$, we consider $\varphi \in X^*$ such that,
$$\varphi(a_k)=0,\, \, \, \, \forall \, \, k \in A.$$
To show that $\varphi$ is equal to $0$ on the whole space $X$, we take  
$$\psi: z\mapsto \varphi(f(z))=\displaystyle{\sum_{j \ge 1}}z^j\varphi(x_j)$$
a holomorphic function such that,
$$\psi(\lambda_k)=\varphi(f(\lambda_k)=\varphi(a_k)=0,\, \, \forall \, \, k \in A.$$
Note that the sequence $(\lambda_j)_{j \ge 1}$ of distinct complex converges to $0$ and is a sequence of zeros for the function $\psi$ which is holomorphic in the disc then by the principle of the isolated zeros $\psi \equiv 0$. Therefore its Taylor coefficients are null and,
$$\varphi(x_j)=0,\, \, \, \forall \, \, j \ge 1.$$
Since the familly $(x_j)_{j \ge 1}$ is dense in $X$ then we have $\varphi \equiv 0$.\\
We now prove the second part of the Lemma, let $T\in \mathcal{L}(X)$ be a bounded linear operator, then
$$Ta_k=Tf(\lambda_k)=T(\displaystyle{\sum_{j \ge 1}}\lambda_k^j x_j)=\displaystyle{\sum_{j \ge 1}}\lambda_k^j Tx_j.$$
If $T =0$, the property is obvious. If $T \ne 0$, there exists $j \ge 1$ such that $Tx_j \ne 0_X$. Let,
$$m=min \{ j:\, Tx_j \ne 0 \}$$
and $Ta_k =\displaystyle{\sum_{j \ge m}}\lambda_k^j Tx_j=\lambda_k^mTx_m+\displaystyle{\sum_{j > m}}\lambda_k^j Tx_j =\lambda_k^m (Tx_m+o(1))$ when $k \rightarrow \infty$ because the sequence $(\lambda_k)_{k \ge 1}$ converges to $0$. The continuity of the norm allows us to have, 
$$\|Ta_k\|=|\lambda_k^m| \big ( \|Tx_m \|+o(1) \big ),\, \, k \rightarrow \infty.$$
Therefore,
$$\frac{Ta_k}{\|Ta_k\|}=\frac{\lambda_k^m (Tx_m+o(1))}{|\lambda_k^m| \big ( \|Tx_m \|+o(1) \big )},\, \, \, k \rightarrow \infty,$$
The result follows. (Note that if $\forall \, \, k,\, \lambda_k>0$, the sequence $\big ( \frac{Ta_k}{\|Ta_k\|} \big )_{k \ge 1}$ has a limit.)
\hfill{$\Box$}
\end{preuve}
This definiton of c.r.c. sequences will be very usefull to prove ours results in the Hardy space $H^2(\mathbb{D},\, X)$ where $X$ is an Hilbert space, and this definition because of the previous remark, contains all the sequences in the case of the finite dimension ($dim\, X< \infty$). Also Lemma \ref{crc} shows that in the case where the dimension of $X$ is infinite, the c.r.c. sequences are not only the sequences who generate finite dimensional subspaces.
\subsection{$S^*$-invariant subspaces generated by lacunary power series and cyclicity}
In this part, we will give a description in $H^2(\mathbb{D},\, X)$ (where $X$ is an Hilbert space) of $S^*$-invariant subspaces generated by lacunary series whose sequence of Taylor coefficients is c.r.c. We first give a necessary condition of cyclicity for any familly of functions $F \subset H^2(\mathbb{D},\, X)$
\begin{lem} \label{ness}
Let $F \subset H^2(\mathbb{D},\, X)$ and $E_F=span(S^{*n}F:\, n \ge 0)$. If $F$ is cyclic then $$X_*(F)=X,$$
where $X_*(F)=\cap_{m \ge 0} span(\widehat{f}(k):\, k \ge m,\, f \in F)$.
\end{lem}
\begin{preuve}
Suppose $X_*(F) \ne X$. Then, there exists $ m \ge 0$ such that  
$$X_m:= span_X(\widehat{f}(k):\, k \ge m,\, f \in F) \ne X.$$
For every $f \in F,\, S^{*m}f \in H^2(\mathbb{D},\, X_m)$ then $\forall \, \, g \in E_F,\, S^{*m}g \in H^2(\mathbb{D},\, X_m)$. Therefore, $\forall \, \, g \in E_F$ we write $g=p+h,\, h\in H^2(\mathbb{D},\, X_m)$ with $deg(p) \le m-1$, and it is clear that $E_F\ne H^2(\mathbb{D},\, X)$.
\hfill{$\Box$}
\end{preuve}
This Lemma presents a necessary condition for cyclicity of a function in $H^2(\mathbb{D},\, X)$. In what follow, we will give the cases where this condition is sufficient for lacunarity. Of course, in the general case under no condition on the spectrum $\sigma(f),\, f \in H^2(\mathbb{D},\, X)$, the assertion $X_*(f)=X$ is not sufficient for cyclicity.\\
\\
\textbf{Notation:}\\
\begin{eqnarray*}
E_f&=&span_{H^2(\mathbb{D},\, X)}(S^{*n} f:\, n \ge 0).\\
X_*(f) &=&\displaystyle \bigcap_{n \ge0} span_X(\widehat{f}(k):\, k \ge n).\\
X_{\infty}(f)&=&span_X \big (\mbox{closed subspaces F of X such that}\, H^2(F)\subset
E_{{S}^{* N}f} \, \, \, \, \forall \, N \ge 0\big
). 
\end{eqnarray*}
In what follows and when there will be no possible ambiguity on the considered function $f$, we write $X_*$ and $X_\infty$ respectively for $X_*(f)$ and $X_\infty(f)$ to simplify the notations. The Theorems proved in this section give us information on the nature of the $S^*$-invariant subspaces generated by lacunary series. Theorem \ref{GTH} shows that these spaces split into two parts, two supplementary subspaces, one of them is a doubly invariant subspace (which means that it is invariant for $S$ and its adjoint $S^*$) and the other one is a finite dimensional $S^*$-invariant subspace generated by a polynomial.
\begin{theo}\label{GTH}
Let $X$ be a separable Hilbert space and $f(z)= \displaystyle{\sum_{k=1}^{\infty}} a_k z^{n_k}\in H^2(\mathbb{D},\, X)$ a lacunary series such that the sequence $(a_k)_{k \ge 1}$ is c.r.c., then
$$E_f=H^2(\mathbb{D},\, X_*) \oplus E_p$$
where $p$ is a polynomial and $p=f-P_{X_*}f$.
\end{theo}
Let $\Lambda \subset \mathbb{N}$ be a sequence of integers. Recall the notation, $$H_\Lambda^2(\mathbb{D},\, X):=\{f \in H^2(\mathbb{D},\, X):\, \sigma(f) \subset \Lambda \}.$$ 
We obtain the following Theorem concerning the cyclicity of c.r.c. lacunary series.
\begin{theo} \label{cyc}
Let $X$ be a separable Hilbert space, $\Lambda \subset \mathbb{N}$ an infinite lacunary set, $F$ a familly of functions in $H_\Lambda^2(\mathbb{D},\, X)$ such that $\forall\, \, f \in F,\, (\widehat{f}(k))_{k \ge 0}$ is a c.r.c. sequence in $X$.\\
The following statements are equivalents. 
\begin{enumerate}[(i)]
\item $span_X(X_*(f):\, f \in F)=X$.
\item $F$ is cyclic in $H^2(\mathbb{D},\, X)$.
\end{enumerate}
\end{theo}
\begin{rem}
The implication $(ii)\Rightarrow(i)$ as we see it before is correct for any familly of functions  $F \subset H^2(\mathbb{D},\, X)$ (see Lemma \ref{ness}).
\end{rem}
If the familly $F$ is made by only one function, $F=\{f\}$,  Theorem \ref{cyc} gives this simple criterion of cyclicity,
\begin{coro} Let $f$ be as in Theorem \ref{cyc},
$$ f \, is \, cyclic \Leftrightarrow X_*(f)=X.$$ 
\end{coro}
Lemmas  \ref{1}-\ref{3} below recall some well known facts on sequences and numerical series that will help us to prove our results.
\begin{lem}\label{1}
If $(n_k)_{k \ge1}$ is a lacunary sequence of non negative integers such that 
$$ n_{k+1} \ge d.n_k \, \, \, \, \, \, \, \, \forall k \ge 1,$$
for some $d>1$, then there exists a number $M$ such that for any integer $N$ it cannot have more than $M$ representations on the form $N=n_j-n_k$.
\end{lem}
\begin{lem}\label{2}
Let $(b_n)_{n \ge 1}$ be a sequence such that for any $n \ge 0,\, b_n \ge 0$ and $\displaystyle{\sum_{n=0}^{\infty}} b_n < \infty$.\\
If $(n_k)_{k \ge0}$ is a lacunary sequence, then 
$$\displaystyle{\sum_{k=1}^{\infty}} \, \displaystyle{\sum_{j>k} b_{n_j-n_k}} < \infty.$$ 
\end{lem}
\begin{lem}\label{3}
If $b_n > 0$, $\displaystyle{\sum_{n=0}^{\infty}} b_n < \infty$ and if  $r_n=\displaystyle{\sum_{k>n}} b_k$, then 
$$\displaystyle{\sum_{n=0}^{\infty}} \frac{b_n}{r_n} = \infty.$$
\end{lem}
First, we need to prove that $X_*=X_{\infty}$.
\begin{lem}\label{4}
Let $X$ be a separable Hilbert space and $f(z)=\displaystyle \sum_{k=1}^{\infty} a_k \, z^{n_k} \in H^2(\mathbb{D},\, X)$ a lacunary series which is not a polynomial and where the sequence $(\frac{a_k}{\|a_k\|})_{k \ge 1}$ is relatively compact.\\
Then, there exists a non-zero element $x \in X$ such that  
$$H^2\otimes x \subset  E_{{S}^{* N}f} \, \, \, \, \forall \, \, N \ge 0.$$
\end{lem}
\begin{preuve}
The proof of this Lemma is an adaptation of the proof given by R. Douglas, H. Shapiro and A. Shields in their paper (see $\cite{DSS}$).\\
Let  $$f(z)= \displaystyle \sum_{k=1}^{\infty} a_k \, z^{n_k}.$$
For any fixed integer $N,\, j \ge 0$, there exists an integer $k_0$ such that for any $k \ge k_0,\, n_k - n_{k-1} \ge j+N$ (which is possible since $f$ is a lacunary series ). We consider
$$\frac{1}{\| a_k \|} S^{*n_k-N-j} S^{*N}f=\frac{a_k}{\| a_k \|}z^j+z^j \displaystyle \sum_{l >k} \frac{a_l}{\| a_k \|} \, z^{n_l-n_k}.$$
Take $r_k=\displaystyle \sum_{l >k} \frac{a_l}{\| a_k \|} \,
z^{n_l-n_k}$. We want to prove that any neighbourhood of $0$
for the weak topologie contains one of the $r_k$. We consider the neighbourhood on the following form
$$\mathit{V}=\big \{ h \in H^2(\mathbb{D},\, X )\, :\, |(h,h_i)|< \,1,\, 1 \le i \le n \big \}.$$
The functions $h_i\, (i=1, \dots n)$ are given elements in $ H^2(\mathbb{D},\, X)$ and $h_i(z)=\displaystyle \sum_{k=1}^{\infty} \hat{h}_i(k)\,  z^{k}$. We obtain 
\begin{eqnarray*}
 |(r_k,h_i)|^2 \le& \bigg(\displaystyle \sum_{l >k} \big| \big(\frac{a_l}{\| a_k \|} \,,\, \hat{h}_i(n_l-n_k)\big) \big| \bigg)^2.\\
\le& \bigg( \displaystyle \sum_{l >k}\frac{\|a_l\|^2}{\| a_k \|^2}\bigg) \bigg(\displaystyle \sum_{l >k}\|\hat{h}_i(n_l-n_k)\|^2 \bigg).\\ 
\end{eqnarray*}
Suppose none of the $r_k$ belongs to $\mathit{V}$, then $1 \le \displaystyle \max_{1 \le i \le n}|(r_k,h_i)| \, \, \, \, \forall \, \, k \ge k_0$. Therefore,
$$\frac{\|a_k\|^2}{\displaystyle \sum_{l >k}\|a_l\|^2} \le \displaystyle \max_{1 \le i \le n} \displaystyle{\sum_{l>k}} \|\hat{h}_i(n_l-n_k)\|^2 \le \displaystyle{ \sum_{i=1}^{n}} \displaystyle{\sum_{l>k}} \|\hat{h}_i(n_l-n_k)\|^2.$$
Taking the sum on $k$, we obtain a contradiction: the left side diverges because of Lemma \ref{3} since the right side converges by Lemma \ref{2}, and so each neighbourhood of $0$
for the weak topologie contains one of the $r_k$. Then there exists a suitable convergent normalized subsequence $(\frac{a_{k_i}}{\|a_{k_i}\|})_{i \ge 1}$ (since $(a_k)_{k \ge 0}$ is c.r.c.), such that the limit $x$ is none zero . The sequence $(r_{k_i})_{i \ge 1}$ does not depend on $N$ neither on $j$, and $x$ verifies
$$xz^j \in E_{S^Nf} \, \, \, \, \forall  \, \, j \ge 0, \, \forall \, \, N \ge 0.$$
Finally,
$$H^2\otimes x \subset  E_{S^{*N}f}\, \, \, \, \forall \, \, N \ge 0.$$
\hfill{$\Box$}
\end{preuve}
Now, we are able to prove the first inclusion  
$$X_{\infty} \subset X_*.$$ 
(see Corollary \ref{incl1} below). Note that this inclusion, Lemma \ref{xy}, Corollaries \ref{coro1} and \ref{incl1} depend on general properties of $S^{*}$-invariant subspaces and not on the lacunarity of a series in $H^2(\mathbb{D},\, X)$ or even on the dimension of $X$.
\begin{lem} \label{xy}
Let $x\, , \,y$ two elements taken in $H^2(\mathbb{D},\, X)$. Then,\\
$$y \in E_x \Rightarrow y_i \in  span(x_j:j \ge 0) \, \, \, \, \forall \, i \ge 0.$$ 
\end{lem}
\begin{preuve}
Let $x,\, y \in H^2(\mathbb{D},\, X)$, then
$$x=\displaystyle \sum_{i=0}^{\infty}x_iz^i\, , \,y=\displaystyle \sum_{i=0}^{\infty}y_iz^i,$$
For any fixed $i_0$,\\
$$S^{*i_0}y=y_{i_0}+\displaystyle \sum_{i>i_0}y_i z^{i-i_0}.$$
Since $y \in E_x$ and this subspace is stable by $S^{*}$ then $ S^{*i_0}y \in
E_x$ therefore there exists a sequence of complex polynomials such that,\\
$$S^{*i_0}y=\displaystyle \lim_{n \to \infty} p_n ( S^{*})x=y_{i_0}+\displaystyle \sum_{i>i_0}y_iz^{i-i_0}.$$
So,
$$y_{i_0}=\displaystyle \lim_{n \to \infty} \big [p_{n} ( S^{*})x\big ](0).$$
Or,
$$[p_n ( S^{*})x \big ](0) \in
span(x_j:\, j \ge 0).$$
Finally, we obtain $y_{i_0} \in  span(x_j:\, j \ge 0\,)$.
\hfill{$\Box$}
\end{preuve}
\begin{coro}\label{coro1} Let $f(z)=\displaystyle {\sum_{j \ge 1}} a_jz^{n_j},\, y=\displaystyle {\sum_{i \ge 1}} y_i z^{i}  \in H^2(\mathbb{D},\, X)$. Then,
$$y \in E_{{S}^{* n_k}f} \Rightarrow y_i \in  span(a_j:j \ge k)\, \,\, \, \forall \, \, i \ge 0.$$
\end{coro}
\begin{preuve}
We take $x={S}^{* n_k}f$ in Lemma \ref{xy}.
\hfill{$\Box$}
\end{preuve}
\begin{coro} \label{incl1} Let $f,\, y$ as in the previous Corollary. Then,
$$y \in H^2(\mathbb{D},\, X_{\infty}) \Rightarrow y_i \in  X_*=\displaystyle \bigcap_{k \ge1}span(a_j:j \ge k)\, \, , \, \forall \, i \ge 0.$$
Therefore,
$$X_{\infty} \subset X_*.$$
\end{coro}
\begin{preuve}
By the definition of $X_{\infty}$, $ H^2(\mathbb{D},\, X_{\infty})\subset \displaystyle \bigcap_{k
  \ge1} E_{{S}^{* n_k}f}$. It suffices to apply Corollary \ref{coro1}.
\hfill{$\Box$}
\end{preuve}
\\
The second step is to prove the inverse inclusion to have the equality $X_*=X_\infty$. For that, we consider the following subspace,
$$A=E_f \ominus H^2(\mathbb{D},\, X_{\infty}).$$ 
And the orthogonal projection $P_A$ on $A$ given by  
$$P_Af=f-P_{H^2(\mathbb{D},\, X_{\infty})}f.$$
\begin{lem}
Let $f(z)=\displaystyle{\sum_{k \ge 0}} a_kz^{n_k}\in H^2(\mathbb{D},\, X)$ a lacunary series, then \begin{enumerate}[(i)]
\item $\sigma (P_A f) \subset \sigma (f)$.
\item $P_Af$ is a polynomial.
\end{enumerate}
\end{lem}
\begin{preuve} Note that, $$H^2(\mathbb{D},\, X_{\infty}):=\big \{ x=\displaystyle \sum_{k=0}^{\infty} x_k z^k \mbox{ such that }x_k \in X_{\infty} \mbox{ and } \displaystyle \sum_{k=0}^{\infty}\|x_k \|^2 < \infty \big \}.$$
In fact, $P_{H^2(\mathbb{D},\, X_{\infty})}$ is a projection defined by componants and $Px:=\displaystyle \sum_{k=0}^{\infty}(P_{X_{\infty}}x_k)z^k$,
where $P_{X_{\infty}}$ is the orthogonal projection on $X_{\infty}$ in $X$.\\
It is clear that if $x\in H^2(\mathbb{D},\, X_{\infty}),\, Px=x$. On the other side, if $x\in H^2(\mathbb{D},\, X_{\infty})^{\bot}=H^2(\mathbb{D},\, X_{\infty}^{\bot})$, then $Px=0$ and,
$$P=P_{H^2(\mathbb{D},\, X_{\infty})}.$$ 
We write,
$$(P_{H^2(\mathbb{D},\, X_{\infty})}f)(z)=\displaystyle \sum_{k=1}^{\infty}(P_{X_{\infty}}a_k)z^{n_k}.$$
which give $(i)$. To prove $(ii)$, suppose $P_Af$ is not a polynomial and denote that $$P_Af \in E_f,$$ 
because $P_Af=f-P_{H^2(\mathbb{D},\, X_{\infty})}f$ and by the definition of $X_\infty$, we have $H^2(\mathbb{D},\, X_{\infty}) \subset E_f$. To simplify the notations, we take
$$y= \displaystyle \sum_{k=0}^{\infty} y_kz^k=P_Af=\displaystyle \sum_{k=1}^{\infty}( P_{(X_{\infty})^{\bot}}a_k)z^{n_k},$$ 
where $y_{n_k}= P_{(X_{\infty})^{\bot}}a_k \in X_{\infty}^{\bot}$.\\
According to what is previous, $P_Af$ is a lacunary series and the sequence $(a_k)_{k \ge 0}$ is c.r.c. then $\big (\frac{y_{n_k}}{\| y_{n_k}\|} \big )_{k \ge 0}$ is relatively compact and since we suppose that it is not a polynomial, $P_Af$ verifies Lemma \ref{4} and,
$$\exists \, \, x \ne 0,\, ( \| x \|=1) \, \, and \, \, x\otimes H^2 \subset span(S^{*k} P_Af:k \ge N), \, \, \, \, \forall \, \, N \ge 0.$$
$x \in {X}_{\infty}^{\bot}$ because $x$ is the limit of a subsequence $ \big (\frac{y_{n_{k_{i}}}}{\| y_{n_{k_{i}}} \|} \big )_{i \ge 1} \in {E}_{\infty}^{\bot} \,$ and
$$x\otimes H^2 \subset \displaystyle \bigcap_{N \ge 0}span(S^{*k} y:k \ge N) \subset \bigcap_{N \ge 0}E_{{S}^{* N}f}.$$
The second inclusion is based on the fact that
$$y=P_Af \in E_f.$$ 
Then ${S}^{* k}y={S}^{* k}P_Af \in E_{{S}^{* k}f} \, \, \, \, \forall \, \, k \ge 0.$ and since $x\otimes H^2 \subset \displaystyle{\bigcap_{N \ge 0}}
E_{{S}^{* N}f}$, we have
$$x\otimes H^2 \subset
H^2(\mathbb{D},\, X_{\infty}),$$ 
Because by definition, $X_{\infty}$ is the maximal subspace such that $H^2(\mathbb{D},\, X_\infty) \subset
E_{{S}^{* Nf}} \, \, \, \, \forall \, \, N \ge 0$. Then,
$$x \in X_{\infty}.$$
Or, we saw before that
$$x \in X_{\infty}^{\bot}.$$ 
Therefore, $x=0$ which is absurd and $P_Af$ is a polynomial.
\hfill{$\Box$}
\end{preuve}
\begin{prop}\label{XX} Let $X$ be a separable Hilbert space, $f \in H_\Lambda^2(\mathbb{D},\, X)$, where $\Lambda$ is a lacunary sequence and $(\widehat{f}(k))_{k \ge 0}$ a c.r.c. sequence. Then,
$$X_* = X_{\infty}.$$
\end{prop}
\begin{preuve}
The inclusion $X_{\infty} \subset X_*$ was proved in Corollary \ref{incl1}. To prove the inverse inclusion, we write $$f=P_Af+P_{H^2(\mathbb{D},\, X_{\infty})}f.$$ 
Let $d$ be the degree of the polynomial $P_Af$, then for any $ k \ge d+1,$
$${S}^{*k}f={S}^{*k}P_{H^2(\mathbb{D},\, X_{\infty})}f \subset H^2(\mathbb{D},\, X_{\infty}).$$
And,
$${S}^{*(d+1)}f=\displaystyle \sum_{k >d} a_k z^{n_k-d-1} \in H^2(\mathbb{D},\, X_{\infty}).$$
Then,
$$ a_k \in X_{\infty} \, \, \, \, \forall \,\,  k>d.$$
Therefore, $$X_* =\displaystyle \bigcap_{n \ge 0} span(a_k:k \ge n) \subset \ span(a_k:k > d) \subset X_{\infty}.$$
\hfill{$\Box$}
\end{preuve}
According to Proposition \ref{XX} and since  
$$X_* = X_{\infty}.$$
We obtain,
$$H^2(\mathbb{D},\, X_*) \subset E_f \subset H^2(\mathbb{D},\, X).$$
\begin{Proof of Theorem 1.1}
If $f$ is a polynomial, the proof is immediate. If $f$ is not a polynomial, we know from what follows that $$f= P_{H^2(\mathbb{D},\, X_{\infty})}f+p,$$ 
where $p$ is a polynomial. By Proposition \ref{XX},
$$X_{\infty}=X_*.$$ 
Recall that by definition, $$H^2(\mathbb{D},\, X_{\infty}) \subset E_f.$$ 
We have to show the double inclusion to prove the equality in the Theorem.\\
It is clear that if $p \in E_f$ then $E_p \subset  E_f$, and on the other hand  $H^2(\mathbb{D},\, X_{\infty}) \subset E_f$ therefore, $$H^2(\mathbb{D},\, X_*) \oplus E_p \subset E_f.$$
Conversely, we have $S^{*n}f=
S^{*n}P_{H^2(\mathbb{D},\, X_{\infty})}f+S^{*n}p \, \, \, \, \forall n
\ge0$ and since $H^2(\mathbb{D},\, X_{\infty})$ is stable by $S^*$,
$$S^{*n}f \subset H^2(\mathbb{D},\, X_*) \oplus E_p \, \, \, \,
\forall n \ge0.$$
And so the inverse inclusion, and the wanted equality $E_f=H^2(\mathbb{D},\, X_*) \oplus E_p$.
\hfill{$\Box$}
\end{Proof of Theorem 1.1}
\begin{rem} Concerning the degree of the polynomial $p$, we observe since $X_{\infty}=X_*$, that there exists a minimal number $N(f)$ such that $$X_{\infty}=span(a_k: \,k \ge N(f)).$$
Moreover, $p= f-P_{H^2(\mathbb{D},\, X_{\infty})}f$ and then, $deg(p)=N(f)-1.$
\end{rem}
\begin{Proof of Theorem 1.2}
Lemma \ref{ness} gives the implication $(ii)\Rightarrow(i)$. To prove that $(i) \Rightarrow (ii)$, we have for any $f \in F$, 
$$E_f=H^2(\mathbb{D},\, X_*(f))\oplus E_p,$$ 
where $p$ is a polynomial. This imply that  
$$H^2(\mathbb{D},\, X_*(f))\subset E_F,$$ 
for any $f \in  F$ and then,
$$E_F \supset span(H^2(\mathbb{D},\, X_*(f)):\,f \in F)=H^2(span(X_*(f):\,f \in F))=H^2(\mathbb{D},\, X).$$
Therefore, $F$ is cyclic.
\hfill{$\Box$}
\end{Proof of Theorem 1.2}

\begin{rem}
To end this part, it is also interesting to see that we can construct by this method cyclic series  in $H^2(\mathbb{D},\, X)$ where $X$ is of infinite dimension. Let $(n_k)_{k \ge 0}$ be a lacunary sequence of non-negative integers, $(e_n)_{n \ge 0}$ be a base of $X$ and take $f(z)=\displaystyle{\sum_{k \ge 1}} a_k z^{n_k}$ where $a_k=(2^{k})^{-\frac{1}{2}}e_{j}$ and $j$ comes from the decomposition of $k$ such that $k=2^{k_0}+j$ for some integers $k_0$ and $j$ with $0\le j \le 2^{k_0+1}-2^{k_0}-1$.\\
Using nearly the same approach as E. Abakumov for the case where the spectrum is a finite union of lacunary sequences (see [\cite{Abak}, p.283]) by taking 
$$g=\frac{\displaystyle{\sum_{1 \le i \le N}} \frac{R_i}{\|a_i\|}S^{n_i}f}{\displaystyle{\sum_{1 \le i \le N}}R_i} \mbox{  with  }  R_k=\frac{\|a_k\|^2}{\displaystyle{\sum_{j >k}}\|a_j\|^2}$$
and since we have for any fixed $j$,
$$\displaystyle{\sum_{k \ge 0}}R_{2^k+j}=\displaystyle{\sum_{k \ge 1}} \frac{\|a_{2^k+j}\|^2}{\displaystyle{\sum_{p >2^k+j}}\|a_p\|^2}=\displaystyle{\sum_{k \ge 0}} \frac{\frac{1}{2^{2^k+j}}}{\displaystyle{\sum_{p>2^k+j}}\frac{1}{2^p}}=\displaystyle{\sum_{k \ge 1}} 1=\infty.$$
Then by following the same steps it is easy to prove that $H^2e_j \in E_f \, \, \forall \,\, j \ge 0$. And $f$ is cyclic.
\end{rem}
\subsection{$S^*$-invariant subspaces generated by a polynomial}
In this part, we want to give a different caracterization of the subspace $E_p$ which appears in Theorem \ref{GTH}.\\ 
\begin{theo}
Let $F \subset H^2(\mathbb{D},\, X)$ a closed subspace. The following statements are equivalent.
\begin{enumerate}[i)]
\item There exists a polynomial $p \in H^2(\mathbb{D},\, X)$ of degree $N$ such that $F=E_p$.\\
\item $dim\, F=N+1$ and $F=K_{\Theta}=H^2(\mathbb{D},\, X)\ominus{\Theta} H^2(\mathbb{D},\, X)$, where $\Theta$ is an inner matricial polynomial of the form $\Theta = \displaystyle{\prod_{k=1}^{\substack{ \curvearrowleft \\ N+1}}} \big( {P_k}^{\perp}+zP_k \big),$ with $dim \, Ker  \, \Theta(0)^{*}=1$ and $P_k: \, X \mapsto X$ are orthogonal projections.
\end{enumerate}
Moreover, the statement $dim\, Ker\, \Theta(0)^*=1$ is equivalent to \\
$\exists \, x_{N+1} \in X$ such that $\|x_{N+1} \|=1$, $P_{N+1}=\langle \, . \, ;x_{N+1}\rangle x_{N+1}$ and $ \, \forall \, \, k=1, \ldots ,N ,$\\
$\exists  \,  x_k \in (1-P_{k+1}) \ldots (1-P_{N+1})X$ such that $\|x_{k}\|=1$ and $P_k=\langle \,. \, ,x_{k}\rangle x_{k}.$
\end{theo} 
\begin{preuve}
Suppose that $F=E_p=span(S^{*n} p\, : \, n \ge  0)$ and $p$ is a polynomial of degree $N$, it is clear that $\displaystyle{ \{ S^{*n} p \} _{i=0}^{N}}$ is a base for $F$ and $dim \, E_p=N+1$. 
$S^*E_p \subset E_p$ and we take $T=S^*\vert_{E_p}$. Since $E_p$ is $S^*$-invariant, it can be represented under the following canonical form $$E_p=H^2(\mathbb{D},\, X) \ominus \Theta H^2(\mathbb{D},\, X),$$
$\Theta$ is an inner matricial function ($z\mapsto \Theta(z)$, $\Theta(\zeta): \, X \mapsto X$ is unitary $|\zeta|=1$).\\
We use the factorization of Blaschke-Potapov,
$$\Theta= V.\, B.\, S,$$
where $V$ is unitary, $S$ is the singular part and $B$ is the following finite Blaschke product,
$$B=\displaystyle{\prod_{k=1}^{\substack{ \curvearrowleft \\ N+1}}} \big( {P_k}^{\perp}+b_{\lambda_k}P_k \big),$$
where $b_{\lambda_k}$ are the usual factors of Blaschke product and $P_k: \, X \mapsto X$ are the orthogonal projections. Moreover, $\sigma (T)=\{\lambda_k,\, 1\le k \le N+1 \}$.\\
In this case, we can have the following reductions; since $dim \, E_p < \infty$ then according to Treil Lemma which is a vectoriel version of Kronecker Theorem (see \cite{Treil}), $\Theta(z)$ is rational. Then the singular part of the factorization is trivial and $S=1$. We have, $\sigma (T)=\{ 0 \}$ because $T^{N+1} \equiv 0$ so $0$ is the only eigenvalue and $b_{\lambda_k}=z$ $\forall \, k$ then $B$ can be written in the simple way, 
$$ B=\displaystyle{\prod_{k=1}^{\substack{ \curvearrowleft \\ N+1}}} \big( {P_k}^{\perp}+zP_k \big).$$
We now prove that $dim \, Ker \,T=1$. We know that $Ker \, T^{N+1}=E_p$ and $Ker \, T^{N} \ne E_p$. If $k \le N$ then $dim \, Ker \, T^{k+1} \ge dim \, Ker \, T^k+1$ or else $Ker \, T^k=Ker T^{k+1}= \ldots=Ker T^{N+1}=E_p $ but $Ker T^N \ne E_p$ and so the contradiction. For $k \le N$, we have $dim \, Ker \, T^{k+1} \ge dim \, Ker \, T+k$ and if $k=N$ then $ dim \,Ker \, T^{N+1}=N+1 \ge  N+dim \, Ker \,T$ which leads to $dim \, Ker \,T \le 1$. But $dim \, Ker \,T \ne 0$ or else $E_p=Ker\, T^{N+1}=\{0\}$ which is impossible. Then $dim \, Ker \, T=1$.\\
In order to show that $dim \, Ker \, \Theta(0)^*=1$, and according to the fact that we can write 
$$Ker(S^*-\overline{\lambda}I)\vert_{K_\theta}=\big \{ \frac{e}{1-\overline{\lambda}z}:\, e \in Ker \, \Theta(\lambda)^*\big \},$$
(see $\cite{Ni1}$), if we take $\lambda=0$, $Ker\, T=Ker \, S^*\vert_{K_{\Theta}}=Ker \, \Theta(0)^*$, then $$dim \, Ker \, T= dim \, Ker \, \Theta(0)^*=1.$$
Conversely, if $\Theta$ is a product of the form $(ii)$ and $T=S^*\vert_{K_\Theta}$ then using the same argument as before $\sigma(T)=\{0\}$. Moreover, using the fact that if any operator $T$ is acting in a finite dimensional space then this operator $T$ is cyclic if and only if $dim \, Ker (T- \lambda I) \le 1 \, \, \, \, \forall \, \, \lambda \in \mathbb{C}$ (see $\cite{NiVa1}$) and it proves that our operator $T$ is cyclic since $Ker\, T=Ker\, \Theta(0)^*$ is of dimension $1$. And there exists $p \in K_\Theta$ such that $$K_\Theta=span(S^{*n}p:\, n \ge 0).$$
Since $\sigma(T)=\{0\}$ and $dim\, K_\Theta=N+1$, it is clear that $T^{N+1}=0$, then $S^{*N+1}p=0$ and $p$ is a polynomial with $deg\, p \le N$. The degree is equal to $N$ because $dim\, K_\Theta=N+1$.\\
We know prove the last statment of the Theorem.\\
Let $\Theta$ be a Blaschke-Potapov product then $\Theta=V.\,B$ and,
$$\Theta(0)=V.\, B(0)=V \displaystyle{\prod_{k=1}^{\substack{ \curvearrowleft \\ N+1}}} \big( 1-P_k\big).$$
Also,
$$\Theta(0)^*=\displaystyle{\prod_{k=1}^{\substack{ \curvearrowright \\ N+1}}} \big( 1-P_k\big).V^*=(1-P_1)\cdots(1-P_{N+1})V^*.$$
So, $Ker \, \Theta(0)^* \supset  Ker \, (1-P_{N+1})V^*$ or $dim \,Ker \, \Theta(0)^*=1$ then $dim \, Ker \, (1-P_{N+1})V^{*}=1$.\\
And since $V$ is a unitary factor then $dim \, Ker \, (1-P_{N+1})=1$ i.e. $Rank \, P_{N+1}=1$. Therefore,
$$dim \, Ker \, \Theta(0)^*=1 \Leftrightarrow Rank \, P_{N+1}=1 \mbox{ and } Ker \,(1-P_1) \ldots (1-P_{N}) \arrowvert_{(1-P_{N+1})X}=\{0 \}.$$
Taking a look at the projections $P_k$ of $\Theta$, we have seen that $Rank \, P_{N+1}=1$ then there exists $x_{N+1} \in X$ such that $\|x_{N+1}\|=1$ and $P_{N+1}=\langle \, . \,;x_{N+1}\rangle x_{N+1}$.\\ 
Since $(1-P_1) \ldots (1-P_{N}) \arrowvert_{(1-P_{N+1})X}$ is injective then $Ker \, (1-P_N) \cap {x_{N+1}}^\perp=\{ 0 \}$ and  $Rank \, P_N \le 1$ and there exists $x_N \in X$ such that  $\|x_{N}\|=1 \, ; \, P_{N}=(.,x_{N})x_{N}$ and $x_{N} \notin {x_{N+1}}^\perp$.\\
By the same way, $(1-P_{N-1})(1-P_{N})\arrowvert_{(1-P_{N+1})X}$ is injective and,
 $$P_{N-1}X \cap (1-P_N)(1-P_{N+1})X=\{ 0 \}.$$
Therefore there exists $x_{N-1} \in X$ such that $\|x_{N-1}\|=1 \, ; \, P_{N-1}=(.,\, x_{N-1})x_{N-1}$ and $x_{N-1} \notin (1-P_N){x_{N+1}}^\perp$. By iteration, we obtain
$$\forall \, \, \, \,  k=1,\ldots N,\, P_k=(.,\, x_{k})x_{k}; \,  \|x_{k}\|=1, \, \, x_k \in (1-P_{k+1})\ldots(1-P_{N+1})X.$$
It is easy to see that this proof is reversible and this property caracterizes the products $\Theta$ participating in $(ii)$.
\hfill{$\Box$}
\end{preuve}
\subsection{Lacunary series and cyclicity for any power of the backward shift $S^*$ in $H^2$}
\begin{lem}\label{equiv}
Let $X$ be a separable Hilbert space, $N \in \mathbb{N}^*$ and,
$$f \in H^2(\mathbb{D},\, X),\, F \in H^2(\mathbb{D},\, X^N),\, f(z)=\displaystyle{\sum_{k \ge 0}} \widehat{f}(k)z^k,\, F(z)=\displaystyle{\sum_{k \ge 0}} \widehat{F}(k)z^k,$$ 
where $X^N=X\times \ldots \times X$ ($N$ times) and $\widehat{F}(k)=\big( \widehat{f}(Nk); \widehat{f}(Nk+1) ; \ldots; \widehat{f}(Nk+N-1) \big) \in X^{N}$.\\
The following statements are equivalent.
\begin{enumerate}[i)]
\item $f$ is $S^{*N}$-cyclic in $H^2(\mathbb{D},\, X)$.\\
\item $F$ is $S^*$-cyclic in $H^2(\mathbb{D},\, X^N)$.
\end{enumerate}
\end{lem}
\begin{preuve}\\
We consider,
\begin{eqnarray*}
\Psi:&H^2(\mathbb{D},\, X) &\longrightarrow H^2(\mathbb{D},\, X^N)\\
&f&\longmapsto \Psi(f)(z)=\displaystyle{\sum_{k \ge0}}\widehat{F}(k) z^k,\, \, avec\, \,  \widehat{F}(k)= \left( \begin{array}{c} \widehat{f}(Nk) \\ \vdots \\ \widehat{f}(Nk+N-1) \end{array} \right).
\end{eqnarray*}
$\Psi$ is an isometric isomorphism and we have the following commutatif diagram,
$$\xymatrix{
    H^2(\mathbb{D},\, X) \ar[r]^{S_X^{*N}} \ar[d]_{\Psi} \ar[rd] & H^2(\mathbb{D},\, X) \ar[d]^{\Psi} \\
    H^2(\mathbb{D},\, X^N) \ar[r]^{S_{X^N}^*} & H^2(\mathbb{D},\, X^N)
  }$$
According to this diagram, $$S_{X^N}^{*}\Psi=\Psi S_X^{*N}.$$ 
Indeed, let $f \in H^2(\mathbb{D},\, X)$ and $g(z)=\displaystyle{\sum_{k \ge N}} \widehat{f}(k)z^{k-N}$. Then,
$$\Psi S_X^{*N}f=\Psi(\displaystyle{\sum_{k \ge N}} \widehat{f}(k)z^{k-N})=\Psi(g)=\displaystyle{\sum_{m \ge 0}} \widehat{G}(m)z^{m}=S_{X^N}^*F(z)=S_{X^N}^*\Psi (f),$$
because $\widehat{G}(m)=\big ( \widehat{g}(Nm),\ldots,\widehat{g}(N(m+1)-1) \big )=\big ( \widehat{f}(N(m+1),\ldots,\widehat{g}(N(m+2)-1) \big ) \, \, \forall \, \, m \ge 0.$\\
Since $S_{X^N}^{*}\Psi=\Psi S_X^{*N}$ then $\forall \, \, k \ge 0, \, , S_{X^N}^{*k}\Psi=\Psi S_X^{*Nk}$ and for any complex polynomial $p$,
$$p(S_{X^N}^*)\Psi(f)=\Psi p(S_X^{*N}).$$ 
If $f$ is cyclic for $S_X^{*N}$ in $H^2(\mathbb{D},\, X)$ then $\Psi(f)$ is cyclic for $S_{X^N}^*$ in $H^2(\mathbb{D},\, X^N)$ and conversely.
\hfill{$\Box$}
\end{preuve}
\\
This connection between cyclicity of the operators $S_{X^N}^*$ and $S_X^{*N}$ gives us a criterion to have lacunary series cyclic for any power of the backward shift operator $S^*$ in $H^2$.
\begin{theo}\label{cycsca}
Let $N \in \mathbb{N}^*$ and $f\in H_\Lambda^2(\mathbb{D},\, X)$, where $\Lambda$ is a lacunary set of non negative integers and $(\widehat{f}(k))_{k \ge 0}$ a c.r.c. sequence in $X$. The following statements are equivalent:
\begin{enumerate}[(i)]
\item $f$ is $S^{*N}$-cyclic in $H^2(\mathbb{D},\, X)$.
\item $span \bigg ( \Big (  \widehat{f}(Nk); \widehat{f}(Nk+1) ; \ldots ; \widehat{f}(Nk+N-1) \Big ): \, \forall \, k \ge m \bigg )= X^{N} \, \, \,\,  \forall \, \, m \ge 0.$
\end{enumerate}
\end{theo}
\begin{preuve}
The function $f$ can be written on the form  
$$f(z)=\displaystyle{\sum_{k \ge 1}} \widehat{f}(k) z^{n_k},$$ 
with $\widehat{f}(k) \ne 0$ and $\frac{n_{k+1}}{n_k} \ge d >1, \, \, \, \, \forall \, \, k \ge 1$.\\
We take for any $k \ge 1$,
$$n_k=Np_k+q_k,$$
where $0 \le q_k \le N-1$. Let\\
$$\widehat{A}(k)=\Big (  \widehat{f}(Np_k); \widehat{f}(Np_k+1) ; \ldots ; \widehat{f}(Np_k+N-1) \Big ).$$
It is easy to verify that  
$$A_k \ne 0 \, \, \, \, \forall \, \, k \ge 1.$$ 
The function $F(z)=\displaystyle{\sum_{k \ge 1}}A_k z^{p_k}$ is a lacunary series $H^2(\mathbb{D},\, X^N)$ because 
$$d<\frac{n_{k+1}}{n_k}=\frac{Np_{k+1}+q_{k+1}}{Np_k+q_k} \le \frac{ Np_{k+1}+N}{Np_k} =\frac{p_{k+1}}{p_k}+\frac{1}{p_k}.$$
Therefore, for any $k_0$ enough bigger to have $\frac{1}{p_k} <\frac{d-1}{2},\, \, k \ge k_0$, we will have $$\frac{p_{k+1}}{p_k}> \frac{d+1}{2}>1,\, \, k \ge k_0.$$ 
And applying Theorem \ref{cyc}, we obtain that $F$ is cyclic if and only if statement $(ii)$ is satisfied. Lemma \ref{equiv} finishes the proof.
\hfill{$\Box$}
\end{preuve}
This previous Lemma is very interesting because we can now construct lacunary series in $H^2$ which are cyclic for any fixed power of the backward shift $S^*$ by giving a necessary and sufficient condition on the Tayor coefficients of the considered function but can also describe the spectrum of $S^{*N}$-cyclic lacunary series in $H^2$.
\begin{coro}
Let $N \in \mathbb{N}^*$ be a fixed non negative integer and $f(z)=\displaystyle{\sum_{k \ge 1}}a_k z^{n_k} \in H^2$ a lacunary series. Suppose that for any $k \ge 0,\, a_k \ne 0$. The following statements are equivalents:
\begin{enumerate}[(i)]
\item $f$ is $S^{*N}$-cyclic.
\item $\forall \, \,m \ge 0,\,  \forall \, \, i=0,\, \ldots,\, N-1, \, \exists \, \, k,\, n_k \ge m : \, n_k \equiv i \, (mod\, N).$
\end{enumerate}
\end{coro}
\begin{preuve}
The criterion on Taylors coefficients in Theorem \ref{cycsca} to have cyclicity for lacunary series is realised because of the nature of the spectrum $\sigma(f)$. Indeed, if $n_k \equiv i\, (mod\, N)$ then $\widehat{F}(k)=a_k e_i$ where $(e_i)_{i =0}^{N-1}$ is a standard base of $\mathbb{C}^N$. 
\hfill{$\Box$}
\end{preuve}
\begin{rem} It is also possible to construct a lacunary series $f(z)=\displaystyle{\sum_{k \ge 1}}a_k z^{n_k} \in H^2$ which is $S^{*N}$-cyclic for any $N \in \mathbb{N}^*$. It suffices for that to consider the lacunary sequence,
$$n_k=(k+1)!+k.$$
For any fixed integer $N>0$, and for any integer $m$, the trunqued sequence $(n_k)_{k \ge m}$ meets all classes modulo $N$ and that guarantees the cyclicity according to Theorem \ref{cyc}.
\end{rem}
\begin{defi}
Let $N$ be a nonzero integer, $cyc \, (S^{*N})$ is the set of cyclic functions in $H^2$ for the operator $S^{*N}$ and for a given function $f \in H^2$, we consider the following set
$$A(f) \stackrel{def}{=} \{N \in \mathbb{N}^*:\, f \in cyc(S^{*N}) \}.$$
\end{defi}
We want to study the nature of $A(f)$ and for that we need a well-known Theorem in number Theory very useful in this situation. It is the Theorem of simultaneous congruences also called the "Chinese Theorem",
\begin{theo}
Let $n_1,\, n_2,\ldots,\, n_k$ be some prime numbers such that each of them is prime with any of the others and $a_1,\, a_2,\ldots,\, a_k$ are any integers, then there exists a unique integer $r$ such that,
$$r \equiv a_i \, (mod\, n_i),\, \, 1 \le i \le k \, \, et  \, \, 0 \le r \le \displaystyle{\prod_{i=1}^{k}} n_i-1.$$
\end{theo}
We give the following caracterization of $A(f)$.
\begin{lem}
Let $A \subset \mathbb{N}^*$ and $1 \in A$, these assertions are equivalents:
\begin{enumerate}[(i)]
\item $N \in A,\, m \mid N \Rightarrow m \in A$.
\item There exists $f \in H^2$ such that $A=A(f)$.
\end{enumerate}
\end{lem}
\begin{preuve}\\
$(ii)\Rightarrow (i):$ It is easy to show that if $f \in cyc \, (S^{*N})$ and if $m \mid N$ then $f \in cyc \, (S^{*m})$.\\
$(i)\Rightarrow (ii)):$ Let $(p_i)_{i \ge 1}$ be the sequence of prime numbers. For every prime number $p$, we consider
$$\alpha_p=max \{\alpha \ge 0:\, \exists\, \, a \in A:\, {p^\alpha} \vert a\},$$
and observe that $0 \le \alpha_p \le \infty$.\\
We take for any $n \ge 1$,
$$a_n=\displaystyle{\prod_{i=1}^n}{p_i}^{\inf(n,\, \alpha_{p_i})},$$
We also consider for any $n \ge 0$, the finite set $B_n$ of prime numbers in $\{p_i\}_{i=1,\ldots n}$ which are not dividing $a_n$. We can apply now the Theorem on simultaneous congruences. We defined the following sequence,
$$n_k=(k+1)!+r_k \, \, \, \, \forall \, \, k \ge 1,$$
where for any $k \ge 1$, $r_k=k \, (mod \, a_k)$ and $r_k=0 \, (mod \, b) \, \, \forall \, b \in B_k$. This sequence is lacunary because  
$$\frac{n_{k+1}}{n_k}=\frac{(k+2)!+r_{k+1}}{(k+1)!+r_k} \ge \frac{(k+2)!}{(k+1)!+(k+1)!} \ge (k+2)\frac{(k+1)!}{2(k+1)!} \ge \frac{k+2}{2},\, \, \, \forall \, \, k \ge 1,$$
and has the particularity that even trunqued from any finite number of it first terms, it takes all classes modulo each number of $A$. Indeed, let $a\in A$, for $n$ enough bigger $a \vert a_n$ since this one contains the decomposition in prime numbers of $a$ and $f$ is $S^{*a_n}$-cyclic because the sequence $r_n$ takes all classes of $a_n$ then $f$ is $S^{*a}$-cyclic since $a$ divides $a_n$.\\
We want to prove that $f$ is not $S^{*a}$-cyclic if $a \not \in A$. Since $a \not \in A$, then according to property $(i)$ in the statement of the Lemma, there exists a prime number $b \mid a$ and $b \in B_n$ if $n$ is enough bigger. The sequence $r_n$ lays on the class $0 \, (mod \, b)$ when $n$ is enough bigger then $f$ is not $S^{*b}$-cyclic and therefore is not $S^{*a}$-cyclique because $b$ divides $a$.
\hfill{$\Box$}
\end{preuve}
\section{Construction of cyclic series whose spectrum is a finite union of lacunary sequences}
In this part, we combinate our approach with the study given by N. K. Nikolski and V. I. Vasyunin $\cite{NiVa2}$ and use some tools they introduce to give a method for constructing cyclic functions $f$ in $H^2(\mathbb{D},\, X)$ with $dim\, X=d < \infty$ and whose spectrum $\sigma(f)=\Lambda$ is a finite union of lacunary sequences. \\
We give here the definition of the degree of cyclicty of a function and some results which are helpfull for us to give the construction we propose here to do.
\begin{defi}
Suppose that $F \subset H^2(\mathbb{D},\, X)$, and $dim\, X=d < \infty$.\\
The degree of cyclicity of the subspace $F$ is defined to be the number 
$$dc(F)=coDim[H^2(\mathbb{D},\, X)\ominus E_F] \stackrel{def}{=}d-\displaystyle{\max_{\zeta \in \mathbb{D}}} \, dim\{g(\zeta):\, g \in H^2(\mathbb{D},\, X)\ominus E_F \}.$$
If $s$ is an integer then let
$$C_s=C_s^X \stackrel{def}{=} \{ F: \, F \subset H^2(\mathbb{D},\, X),\, dc(F) \le s \}.$$ 
\end{defi}
\begin{coro} \label{221} $\cite{NiVa2}$ 
Let $F$ and $G$ be subspaces of $H^2(\mathbb{D},\, X)$.
\begin{enumerate}[i)]
\item If $F \subset G$, then $dc(F) \le dc(G)$.
\item $dc(F+G) \le dc(F)+dc(G)$.
\item If $E_F=H^2(\mathbb{D},\, X) \ominus \Theta_F H^2(\mathbb{D},\, X)$ is the canonical representation of the space $E_F$, then 
$$dc(F+G)=dc(F)+dc(P_+ \Theta_F^*G).$$
\end{enumerate}
\end{coro}
\begin{theo} \label{312} $\cite{NiVa2}$ 
Let $\psi = \big \{ \psi_i \big \}_{1 \le i \le d} \in H^2( \mathbb{C}^d)$. The following assertions are equivalent:
\begin{enumerate}[i)]
\item $\psi$ is a cyclic vector in $H^2(\mathbb{D}, \mathbb{C}^d).$
\item For any cyclic vector $\phi$ in $H^2( \mathbb{D},\, \mathbb{C}^{d-1})$, there exists an index $i$ $( 1 \le i \le d)$ such that the vector $\binom{\phi}{\psi_i}$ is cyclic in $H^2( \mathbb{C}^d).$ 
\end{enumerate}
\end{theo}
In the scalar case $H^2$, as we mentioned it before, E. Abakumov proved that functions whose spectrum is a finite union of lacunary sequences are cyclic. We will give some constructions of cyclic vector-valued functions whose spectrum is a finite union of lacunary sequences. \\
Here, we give a Theorem on the existence of cyclic functions for the backward shift in $H^2(\mathbb{D},\, X)$ with $dim \, X < \infty$ and whose coordinate functions are lacunary series in $H^2$ and with the possibility to choose this lacunarity for each one of them.
\begin{theo}
Let $\Lambda_j \subset \mathbb{N},\, i=1,\ldots ,\, d$ arbitrary lacunary sequences. There exists a cyclic function $\psi = \big \{ \psi_i \big \}_{1 \le i \le d} \in H^2(\mathbb{D},\, \mathbb{C}^d)$ such that $\sigma(\psi_i)=\Lambda_i,\, 1 \le i \le d$.
\end{theo}
\begin{preuve}
We consider $\forall \, \, i, \,1 \le i \le d$ any infinite lacunary sets $\Lambda_i$ of $\mathbb{N}$.\\
We take $F_1=\psi_1$ and $\sigma(\psi_1)=\Lambda_1$ which is cyclic in $H^2$ according to the Theorem of Douglas-Shapiro-Shields.\\
Let $\Psi_2=\binom{\psi_{2,\, 1}}{\psi_{2,\, 2}}$ where $\sigma(\psi_{2,\,i})= \Lambda_2,\, i=1,\, 2$. $\Psi_2$ is cyclic if the criterion of Theorem \ref{cyc} is verify i.e.  
$$span\big( \binom{\widehat{\psi}_{2,\, 1}(k)}{\widehat{\psi}_{2,2}(k)} \, : \, k \ge N \big)= \mathbb{C}^2,\, \, \, \, \forall \, \, N \ge 0.$$
This condition is easily realised for some convenient coefficients. Since both of the functions $F_1$ and $\Psi_2$ are cyclic we can apply Theorem \ref{312}. Therefore, there exists $\psi_2= \psi_{2,\, j}$ for some index $j \in \{1;\, 2 \}$ such that the function $F_2=\binom{\psi_1}{\psi_2}$ is cyclic in $H^2(\mathbb{D},\, \mathbb{C}^2)$. \\
In the same way, we can consider the previous cyclic function $F_2$ and 
$$\Psi_3= \left( \begin{array}{c} \psi_{3,\, 1} \\ \psi_{3,\, 2} \\ \psi_{3,\, 3} \end{array} \right),$$ where $\forall \, \, \, \, i\in \{1,\, 2,\, 3 \}, \, \sigma(\psi_{3,\, i})=\Lambda_3$ is such that,
$$span\Big( \left( \begin{array}{c} \widehat{\psi}_{3,\,1}(k) \\ \widehat{\psi}_{3,\, 2}(k) \\ \widehat{\psi}_{3,\,3}(k) \end{array} \right)  \, : \, k \ge N \big)= \mathbb{C}^3 \, \, \, \, \forall \, \, N \ge 0.$$
This last condition is easily realised with appropriate coefficients. We apply again Theorem \ref{312} which gives the existence of $\psi_3= \psi_{3,\,j}$ for some index $j \in \{1,\, 2,\, 3 \}$ such that the function $$F_3= \left( \begin{array}{c} \psi_1 \\ \psi_2 \\ \psi_3 \end{array} \right)$$ is cyclic in $H^2(\mathbb{D},\, \mathbb{C}^3)$ with $\sigma(\psi_i)=\Lambda_i,\, 1\le i \le 3$. By iteration, we obtain a cyclic function $\Psi= \big \{ \psi_i \big \}_{1 \le i \le d}$ in $H^2( \mathbb{C}^d)$ and its coordinate functions $\psi_i$ are lacunary series in $H^2$ wiht spectrum in $\Lambda_i$.
\hfill{$\Box$}
\end{preuve}
We want in this part to extend Theorem \ref{cyc} for some particular cases of functions in $H^2(\mathbb{D},\, X^d)$  whose coordinate functions are lacunary series with different spectrum. We start by giving a general Lemma without any condition on the spectrum of the given functions.
\begin{lem} \label{lemcyc}
Let $X$ be a separable Hilbert space and $f \in H^2( \mathbb{D},\, X^d)$,\\
The following statements are equivalent:
\begin{enumerate}[i)]
\item $f= \left( \begin{array}{c} f_1 \\ f_2 \\ \vdots \\ f_d \end{array} \right)$  is cyclic. 
\item $\forall \, \vartheta_i  \in H^{\infty} ,\, \vartheta_i \not \equiv 0,  \, \, \,i=1,\ldots d, \,  \, \, 
\left( \begin{array}{c} P_+(\overline{\vartheta}_1 f_1) \\ P_+(\overline{\vartheta}_2 f_2) \\ \vdots \\ P_+(\overline{\vartheta}_d f_d) \end{array} \right)$ is cyclic.
\item $\exists \, \vartheta_i  \in H^{\infty} ,\, \vartheta_i \not \equiv 0,  \, \, \,i=1,\ldots d, \,  \, \, 
\left( \begin{array}{c} P_+(\overline{\vartheta}_1 f_1) \\ P_+(\overline{\vartheta}_2 f_2) \\ \vdots \\ P_+(\overline{\vartheta}_d f_d) \end{array} \right)$ is cyclic.
\end{enumerate}
\end{lem}
\begin{preuve} Since $ii)\Rightarrow iii)$ is trivial, we first prove that $i)\Rightarrow ii)$ and then $iii)\Rightarrow i)$.\\
$i)\Rightarrow ii)$. $\forall \, \, \vartheta_i  \in H^{\infty} \, \vartheta_i \not \equiv 0,  \, \, \,i=1,\ldots d $, we define the operator $\mathcal{D}$ by,
\begin{eqnarray*}
\mathcal{D}:  H^2(\mathbb{D},\, X^d) &\longrightarrow& H^2(\mathbb{D},\, X^d) \\ 
g &\longmapsto& \mathcal{D}(g)=\left( \begin{array}{c} P_+(\overline{\vartheta}_1 g_1) \\ P_+(\overline{\vartheta}_2 f_2) \\ \vdots \\ P_+(\overline{\vartheta}_d g_d) \end{array} \right).
\end{eqnarray*}
We have to show that $\mathcal{D}H^2(\mathbb{D},\, X^d)$ is dense in $H^2(\mathbb{D},\, X^d)$ and it suffices to prove that 
$$Ker \,\mathcal{D}^*=\{ 0 \}.$$
To define $\mathcal{D}^*$, we look at the scalar product,
$$\langle \mathcal{D}h;g\rangle=\displaystyle{\sum_{i=1}^d}\langle P_+(\overline{\vartheta}_i h_i);g_i\rangle=\displaystyle{\sum_{i=1}^d}\langle h_i;\vartheta_ig_i\rangle=\langle h;\displaystyle{\sum_{i=1}^d}\vartheta_ig_i\rangle=\langle h;\mathcal{D}^*g\rangle  \, \, \, \forall g, \, h \in H^2(\mathbb{D},\, X^d).$$
We obtain, $$\mathcal{D}^*g=\left( \begin{array}{c} \vartheta_1 g_1 \\ \vartheta_2 g_2 \\ \vdots \\ \vartheta_d g_d \end{array} \right)\, \, \, \,   \forall \, \, g \in H^2(\mathbb{D},\, X^d).$$
If $\mathcal{D}^*g=0$, then $\vartheta_i g_i=0 \, \, \, \, \forall \, i=1,\ldots d$. But $\vartheta_i \ne 0  \, \, \, \,  \forall \, \, i=1,\ldots d$.\\
Therefore, $g_1=\ldots=g_d=0$ and $g=\left( \begin{array}{c} g_1 \\ g_2 \\ \vdots \\ g_d \end{array} \right)=0.$\\
And so $\mathcal{D}$ is dense in $H^2(\mathbb{D},\, X^d)$. Moreover, for every $n \ge 0$,
$$S^{*n}\mathcal{D}f=\mathcal{D}S^{*n}f.$$
If $f$ is a cyclic function in $H^2(\mathbb{D},\, X^d)$, then
$$E_{\mathcal{D}f}=span(S^{*n}\mathcal{D}f:\, n \ge 0) \supset \mathcal{D}(span(S^{*n}f:\, n \ge 0))=\mathcal{D}E_f.$$
But $E_f=H^2(\mathbb{D},\, X^d)$ and $E_{\mathcal{D}f}$ is closed then, $$E_{\mathcal{D}f}=H^2(\mathbb{D},\, X^d).$$
We verify that $iii) \Rightarrow i)$. Let $\vartheta_i  \in H^{\infty},\, \vartheta_i \ne  0$ and $f \in H^2(\mathbb{D},\, X^d)$ such that $\mathcal{D}f$ is cyclic.\\
We take, $$\theta=\vartheta_1\vartheta_2\ldots\vartheta_d.$$
And for every $i=1,\ldots d$,
$$\theta_i=\vartheta_1 \ldots \stackrel{\vee}{\vartheta_i} \ldots \vartheta_d.$$
Suppose that, $\left( \begin{array}{c}P_+(\overline{\vartheta}_1f_1) \\ P_+ (\overline{\vartheta}_2f_2) \\ \vdots \\ P_+ (\overline{\vartheta}_df_d) \end{array} \right)$ is cyclic.\\
If we consider the previous application $\mathcal{D}$ defined with the functions $\theta_1,\ldots,\, \theta_d$ and using the already proved implication $i)\Rightarrow ii)$, then the function 
$$\mathcal{D}\left( \begin{array}{c}P_+(\overline{\vartheta}_1f_1) \\ P_+ (\overline{\vartheta}_2f_2) \\ \vdots \\ P_+(\overline{\vartheta}_df_d) \end{array} \right)=P_+ \overline{\theta} \left( \begin{array}{c}f_1 \\ f_2 \\  \vdots \\  f_d \end{array} \right)=P_+\overline{\theta}f$$ is cyclic. By approximating $\theta$ by its Fej\'{e}r polynomials $\phi_n$, we obtain $\|\phi_n\|_\infty \le \|\theta\|_\infty$ and $\phi_n(\zeta)\rightarrow \theta(\zeta)$ $a.e. \, \zeta \in \mathbb{T}$ which leads to 
$$\displaystyle{\lim_{n \rightarrow + \infty}} \| \phi_n(S^*)f-P_+\overline{\theta}f\|_2=0.$$
Then $P_+\overline{\theta}f \in E_f$ and so $E_f=H^2(\mathbb{D},\, X^d)$.
\hfill{$\Box$} 
\end{preuve}
\begin{theo} \label{zed}
Let $X$ be a separable Hilbert space and $\Lambda$ an infinite lacunary set.
Choose $ \{m_1,\, m_2, \dots, \, m_d \} \subset \mathbb{Z}$ some fixed integers.\\
And let, $$f = (f_i )_{1 \le i \le d} \in H^2(\mathbb{D},\, X^d),$$ 
where $\sigma(f_i) \subset
\Lambda+m_i$. Suppose that the sequence
$$(\widehat{f}_1(k+m_1),\,\widehat{f}_2(k+m_2),\ldots,\,\widehat{f}_d(k+m_d)))_{k \ge 0}$$ is c.r.c. in $X^d$, then
$$f \, \, cyclic \Leftrightarrow span\Big( \left( \begin{array}{c}
    \widehat{f}_1(k+m_1) \\ \widehat{f}_2(k+m_2) \\ \vdots \\ \widehat{f}_n(k+m_d) \end{array} \right): \, k \ge N \big)= X^d \, \, \, \, \forall \, \, N \ge 0.$$
\end{theo}
\begin{preuve}
We use the criterion of cyclicity from Theorem \ref{cyc} for lacunary series in $H^2(\mathbb{D},\, X^d)$. We suppose, without loss of generality that $\sigma(f_d)=\Lambda=\{ n_k
\}_{k \ge 1}$ et $m_1 \ge m_2 \ge \dots \ge m_d=0$.\\
It is a direct application of Lemma \ref{lemcyc} which gives the equivalence. Take  $\vartheta_i=z^{m_i}, \, i=1,\ldots d$. Then,
$$\mathcal{D}f=\left( \begin{array}{c} S^{*m_1}f_1 \\ S^{*m_2}f_2 \\ \vdots \\ S^{*m_d}f_d \end{array} \right)$$ 
is a lacunary series because the spectrum of each coordinate function is included in $\Lambda$. Moreover, $$\widehat{\mathcal{D}f}(k)=(\widehat{f}_i(k+m_i))_{i=1}^{d},\, \, \, \, \forall \, \, k \ge 0.$$
The application of Theorem \ref{cyc} and Lemma \ref{lemcyc} finish the proof.
\hfill{$\Box$}
\end{preuve}
\textbf{Example:}
Let $(n_k)_{k \ge 1}$ be a lacunary sequence of Hadamard and $f \in H^2(\mathbb{D},\, \mathbb{C}^2)$ such that,
$$f(z)= \displaystyle{\sum_{k \ge 1}} \left( \begin{array}{c} a_k \\ 0  \end{array} \right) z^{n_k}+\displaystyle{\sum_{k \ge 1}} \left( \begin{array}{c} 0 \\ b_k  \end{array} \right) z^{n_k+1}.$$ 
According to the previous Proposition, if $\forall \, \, m \ge 1$, 
$$span( \left( \begin{array}{c} a_k \\ b_k  \end{array} \right):\, k\ge m)=\mathbb{C}^2,$$
then $f$ is cyclic.
\\
\\
\\
It is well known among properties of cyclic vectors in $H^2$ the following fact
 $$f\in C,\, g \in N \Rightarrow f+g \in C.\, \, \, \, (*)$$
In general, this is not true in the space $H^2(\mathbb{D},\, X)$ where $X$ is a separable Hilbert space as we can see in the following trivial example: If we consider in $H^2(\mathbb{D},\, \mathbb{C}^2)$ a cyclic function of the form $$f=\left( \begin{array}{c} f_1 \\ f_2 \end{array} \right),$$
and the function $g$, 
$$g=\left( \begin{array}{c} -f_1 \\ 0 \end{array} \right).$$
Then, $dc(f)=2,\, dc(g)=1$ and the sum,
$$f+g=\left( \begin{array}{c} 0 \\ f_2 \end{array} \right)$$ 
is not cyclic and so $f\in C,\, g \in N \nRightarrow f+g \in C.$ \hfill{$\Box$}\\
\\
In the case of vector-valued spaces $H^2(\mathbb{D},\, \mathbb{C}^d)$, and using the degree of cyclicity of a function $f \in H^2(\mathbb{D},\, \mathbb{C}^n)$, it is possible to save a part of the previous implication $(*)$. For example if $f,\, g \in H^2(\mathbb{D},\, \mathbb{C}^d)$ are such that $f \in C $ (then $dc(f)=d$) and $g \in N$ with $dc(g)=0$ then
$$f+g \in C.$$
We can see this fact as a direct consequence of $(iii)$ in the Corollary \ref{221}.\\
We will describe some pair of functions $\{f,\, g\} \subset H^2(\mathbb{D},\, \mathbb{C}^d)$, $f$ cyclic and $g$ non-cyclic with $dc(g)=1$, such that $f+g$ is cyclic.
The following Proposition allows us to construct cyclic functions in $H^2(\mathbb{D},\, \mathbb{C}^2)$ whose spectrum is a finite union of the form $(n_k)_{k \ge 1}\cup (n_k+1)_{k \ge 1}$ where $(n_k)_{k \ge 1}$ is  lacunary sequence.
\begin{prop}
Let $(n_k)_{k \ge 1}$ be a lacunary sequence of Hadamard.\\
$\varphi_1= \left( \begin{array}{c} f_1 \\ g_1 \end{array} \right),\, \varphi_2=\left( \begin{array}{c} f_2 \\ g_2 \end{array} \right)$ two functions in $H^2(\mathbb{D},\, \mathbb{C}^2)$ such that $\sigma(f_1),\, \sigma(g_1) \subset (n_k)_{k \ge 1}$ and $\sigma(f_2),\, \sigma(g_2) \subset (n_k+1)_{k \ge 1}$. Suppose that $dc(\varphi_1)=2$ and $dc(\varphi_2)=1$ and the following condition fullfiled 
$$span(\left( \begin{array}{c} \widehat{f}_1(k) \\ \widehat{g}_1(k) \\ \widehat{g}_2(k-1) \end{array} \right):\, k \ge m )=\mathbb{C}^3,\, \, \, \, \forall\, \, m \ge 1.$$
Then $\varphi_1+\varphi_2$ is cyclic.
\end{prop}
We omit the proof here because it is nearly the same as the following one we give for the next Proposition \ref{propsuffi} with minors changements.\\
\\
\textbf{Remark:} We can generalize this Proposition in $H^2(\mathbb{D},\, \mathbb{C}^d)$ to series of the form $f=\displaystyle{\sum_{k=1}^r} \varphi_k$ where $\sigma(\varphi_i) \subset (n_k+i)_{k \ge 1},\, i=1,\ldots r$ and $dc(\varphi_1)=d,\, dc(\varphi_i)=1$ for $i=2,\ldots r$ for any $d,\, r \ge 1$. Note that if $d=1$, the Proposition is true because of the result proved by E. Abakumov on the cyclicicty of functions whose spectrum is a finite union of lacunary sequences.\\
\vspace{3mm}
We can see by the following example that the criterion of Theorem \ref{cyc} does not work for functions in $H^2(\mathbb{D},\, \mathbb{C}^d)$ whose spectrum is finite union of lacunary sequences with bounded blocks. The next Proposition (as Theorem \ref{zed}) is a partial answer by giving a sufficient condition of cyclicity for these functions.\\
\\
\textbf{Example:}
$$F=\left( \begin{array}{c} f \\ S^*f \end{array} \right)\in H^2(\mathbb{D},\, \mathbb{C}^2).$$
The spectrum of this function is the union of two lacunary sequences with bounded blocks of length two and it is clear that the criterion of Theorem \ref{cyc} is realised because for every $m \ge 0$,
$$span(\widehat{F}(k):\, k \ge m)=span \big(\left( \begin{array}{c} \widehat{f}(k) \\ 0 \end{array} \right);\, \left( \begin{array}{c} 0 \\ \widehat{f}(k) \end{array} \right):\, k\ge m \big )=\mathbb{C}^2.$$
But it is easy to see that $F$ is not cyclic. 
\begin{prop}\label{propsuffi}
Let $r,\, d \ge 1$ two fixed integers and $\varphi=\displaystyle{\sum_{k=1}^r} \varphi_k \in H^2(\mathbb{D},\, \mathbb{C}^d),$ \\
where $\sigma(\varphi_i) \subset (n_k+i-1)_{k \ge 1},\, i=1,\ldots r$. Suppose that for every $m \ge 0$,
$$span(\left( \begin{array}{c} \widehat{\varphi}_1(k)\\ \widehat{\varphi}_2(k+1) \\ \vdots \\ \widehat{\varphi}_r(k+r-1) \end{array} \right):\, k \ge m )=\mathbb{C}^{d\times r}.$$
Then $\varphi$ is cyclic.
\end{prop}
\begin{preuve}
We take in $H^2(\mathbb{D},\, \mathbb{C}^{d\times r})$ the following subspaces,
$$F= span \left( \begin{array}{c} \varphi_1+\varphi_2+\varphi_3+\ldots+ \varphi_r\\ 0 \\ \vdots \\0 \end{array} \right) ,\, G= span \Bigg(\left( \begin{array}{c} \varphi_2 \\ -\varphi_2 \\ 0 \\ \vdots \\ 0 \end{array} \right);\, \left( \begin{array}{c} \varphi_3 \\ 0 \\ -\varphi_3 \\ 0 \\ \vdots \\0 \end{array} \right);\ldots ;\, \left( \begin{array}{c} \varphi_r \\ 0 \\ \vdots \\ 0 \\ -\varphi_r \end{array} \right) \Bigg).$$
Using the equation of Corollary \ref{221} and the fact that $dc(P_+ \Theta_F^*G) \le dc(G)$ (it is a consequence of DC Lemma in $\cite{NiVa2}$) and because $G$ is generated by $r-1$ functions of degree $d$, we obtain $$dc(G) \le (r-1)d.$$
Moreover, if we suppose that $dc(F)<d$ then $dc(F+G)<rd$. On the other hand, 
$$\left( \begin{array}{c} \varphi_1 \\ \varphi_2 \\ \vdots \\ \varphi_d \end{array} \right) \in F+G.$$ 
This function is cyclic in $H^2(\mathbb{D},\, \mathbb{C}^{d\times r})$ according to the condition given in the Proposition and because of Theorem \ref{zed}, then $dc(F+G)=d \times r$. This contradiction shows that in fact $\varphi$ is cyclic.
\hfill{$\Box$}
\end{preuve}
\section{Lacunary series with bounded blocks}
We want here to find a criterion of cyclicity which is also available for some type of functions who generalize lacunary series. The following Lemma and Proposition can be regard as preparatories for such a description.
\begin{defi}
Let $X$ be a separable Hilbert space, $f \in H^2(\mathbb{D},\, X)$ has the bounded block lacunary property if there exists $N \in \mathbb{N}$ such that 
$$\Psi_N(f)\stackrel{def}{=}\displaystyle{\sum_{k \ge0}}\widehat{\Psi}_N(f)(k) z^k,$$
where $\widehat{\Psi}_N(f)(k)= \left( \begin{array}{c} \widehat{f}(Nk) \\ \vdots \\ \widehat{f}(Nk+N-1) \end{array} \right) \in X^N$ and $\sigma(\Psi_N(f))$ is a lacunary sequence.
\end{defi}
\textbf{Examples :}
\begin{enumerate}[(i)]
\item The sequence $(a^k)_{k \ge 0}$ where $a \in \mathbb{N}\setminus \{0,\, 1\}$ has the lacunary bounded block property.
\item The sequence $m_k=k!+k$ has not the lacunary bounded blocks property. This example also shows that this property is not connected to the speed of growth of the sequence.
\item The sequence $\cup_{k \ge 1}[n_k,\, n_k+N]$ where $(n_k)_{k \ge 1}$ is a lacunary sequence has not the property in general.
\end{enumerate}
We give a criterion of cyclicity for series with the property of lacunary bounded blocks.
\begin{lem}\label{01..N}
Let $F$ be a familly of functions in $H^2(\mathbb{D},\, X)$. The following statements are equivalents. 
\begin{enumerate}
 \item $F$ is cyclic in $H^2(\mathbb{D},\, X)$.
 \item The familly $\{ \Psi_N(F),\, \Psi_N(S^*F),\ldots,\, \Psi_N(S^{*N-1}F) \}$ is cyclic in $H^2(\mathbb{D},\, X^N)$.
\end{enumerate}
\end{lem}
\begin{preuve}
It suffices to see that from one hand $\Psi_N$ is a isometric ismorphism from $H^2(\mathbb{D},\, X)$ to $H^2(\mathbb{D},\, X^N)$ and on the other hand any integer $k$ can be writen with the euclidian division $k=p_kN+q_k$, where $0 \le q_k \le N-1$ and therefore,
$$\Psi_N(S_X^{*k}F)=S_{X^N}^{*p_k}\Psi_N(S^{*q_k}f).$$
\hfill{$\Box$}
\end{preuve}
\begin{prop}
Let $N$ be an integer and $F \in H^2(\mathbb{D},\, X)$ a finite familly of functions such that $\forall \, \, f \in F$, the sequence $(\widehat{f}(k)_{k \ge 0}$ is c.r.c. and $\Psi_N(f)$ has the lacunary bounded blocks property. The following assertions are equivalent:
\begin{enumerate}
 \item $F$ is cyclic in $H^2(\mathbb{D},\, X)$.
 \item $X_*(\Phi)=X^N$ where $\Phi=\{\Psi_N(F),\, \Psi_N(S^*F),\ldots,\,\Psi_N(S^{*N-1}F)\}$.
\end{enumerate}
\end{prop}
\begin{preuve}
It is a consequence of Lemma \ref{01..N} and Theorem \ref{cyc}.
\hfill{$\Box$}
\end{preuve}
In this part, we want to solve the problem of cyclicity for series whose spectrum is included in infinite sets of the form 
$$\Lambda = \displaystyle{\cup_{k \ge 0}} [n_k; n_k+N],$$ 
where $(n_k)_{k \ge 0}$ is a lacunary sequence and $N$ a fixed integer.\\
We take such a series $f$ and we give the construction of an $S$-invariant subspace $F \subset E_f$ which play the role of the subspace $H^2(\mathbb{D},\, X_*)$ in Theorem \ref{GTH}. At this time, we can not prove that $F$ is $S^*$-invariant (and in general it is not the case), but we still think that the final criterion depends directly on the nature of $F$. Here is the construction.
\begin{lem}\label{bloc}
Let $X$ be a separable Hilbert space and
$$f(z)= \displaystyle{\sum_{k \ge 1}} P_k(z)z^{n_k} \in H^2(\mathbb{D},\, X),$$ 
where $P_k$ are polynomials of degree less or equal to $N$ and $(n_k)_{k \ge 1}$ is a lacunary sequence. Suppose that the sequence $\big(\frac{P_k}{\|P_k\|} \big)_{k \ge 1}$ is relatively compact in $\mathcal{P}_N(X)$.\\
Then, there exists a non zero polynomial $P,\, deg\, P\le N$ with values in $X$ such that 
$$ H^2\otimes P \subset E_{S^{*n}f} \, \, \, \, \forall \, \, n \ge0.$$
\end{lem}
We omit the proof which is closed to the one we give in Lemma \ref{4} and it is clear from it that $deg(P) \le N$. We call $\mathcal{P}_N(X)$, the set of polynomials whose degree is less or equal to $N$ and with values in $X$.
\begin{coro}\label{corobloc}
Let $f$ as in Lemma \ref{bloc}. Then there exists a closed subspace $L\subset \mathcal{P}_N(X)$ such that 
\begin{enumerate}[(a)]
\item $L\ne \{0\}$.
\item $F:=H^2\otimes L \subset \displaystyle{\bigcap_{k \ge 0}}E_{S^{*k}f}\subset E_f$.
\item If $g \in \mathcal{P}_N(X)$ et $H^2\otimes g \subset E_f$, then $g \in L$ (``$L$ is maximal'').
\end{enumerate}
\end{coro}
\begin{preuve}
Indeed, by definition, $H^2\otimes L=clos_{H^2(\mathbb{D},\, X)} \big \{\displaystyle{\sum} h_jp_j:\, h_j \in H^2,\, p_j \in L\big \}.$\\
Then, if $L_\alpha,\, \alpha \in A$, satisfied the statement $(b)$, then it is true for $span_{H^2(\mathbb{D},\, X)}(L_\alpha:\, \alpha \in A)$. The rest are direct consequences of Lemma \ref{bloc}.
 \hfill{$\Box$} 
\end{preuve}
The following Theorem gives the others properties of the subspace $F=H^2\otimes L$ defined in Corollary \ref{corobloc} which are connected to the possible cyclicity of $f$.
\begin{theo}\label{teobloc}
Let $X$ be a separable Hilbert space and
$$f(z)= \displaystyle{\sum_{k \ge 1}} P_k(z)z^{n_k} \in H^2(\mathbb{D},\, X),$$ 
where $P_k$ are the polynomials with degree less or equal to $N$ and $(n_k)_{k \ge 1}$ a lacunary sequence. Suppose that the sequence $(P_k)_{k \ge 1}$ is c.r.c. in $\mathcal{P}_N(X)$. Let $F=H^2\otimes L$, the subspace in Corollary \ref{corobloc}. Then,
\begin{enumerate}[(i)]
\item $SF \subset F \subset \displaystyle{\bigcap_{k \ge 0}}E_{S^{*k}f}\subset E_f$.
\item $f=g+p$ where $g\in F$ and $p$ is a polynomial.
\item $\displaystyle{\bigcap_{k \ge 0}}E_{S^{*k}f}=clos(F+S^*L+\ldots+S^{*N}L)$ and there exists $d \ge 0$ such that $$S^{*d}f \in clos(F+S^*L+\ldots+S^{*N}L)=\displaystyle{\bigcap_{k \ge 0}}E_{S^{*k}f}.$$
In particular,
$$S^{*d}E_f=clos(F+S^*L+\ldots+S^{*N}L).$$
\item If $dim\, X<\infty$ then the necessary and sufficient condition to have $f$ cyclic is for $L$ to have the maximal rank in $X$, i.e.
$$dim\, X=r:=\displaystyle{\max_{|z|<1}}\, dim(p(z):\, p\in L).$$
\end{enumerate}
 \end{theo}
\begin{preuve}
It is clear that we have $(i)$.\\
$(ii)$ By definition, $f=\displaystyle{\sum_{k \ge 1}}z^{n_k}P_k$, with $deg(P_k) \le N$. On the other hand, every convergent sum on the form $g=\displaystyle{\sum_{k \ge 1}}z^{n_k}Q_k$, where $Q_k  \in L$ is in $F=H^2\otimes L$. In particular, it is the case for the series $g$ with $Q_k=P_L(P_k)$. Then, $p:=f-g \in E_f$ and $p=\displaystyle{\sum_{k \ge 1}}z^{n_k}P_{L^\perp}(P_k)$ where $L^\perp=\mathcal{P}_N(X)\ominus L$. \\
Lemma \ref{bloc} can be apply to the series $p$ instead of 
$f$. If $p$ is not a polynomial, we obtain a polynomial $R\ne 0,\, R \in \mathcal{P}_N(X)\ominus L$ such that $H^2\otimes R \subset \displaystyle{\bigcap_{k \ge 0}}E_{S^{*k}f}\subset E_f$. The contradiction (with Corollary \ref{corobloc}) shows $p$ is a polynomial. By taking $g=\displaystyle{\sum_{k \ge 1}}z^{n_k}P_L(P_k) \in H^2 \otimes L$ and $f=g+p$, we obtain the result.\\
$(iii)$ By definition,
$$H^2\otimes L=clos_{H^2(\mathbb{D},\, X)}\big \{\displaystyle{\sum} c_jz^jp_j:\, p_j \in L\big \},$$
(The sums are finite). On the other hand, for every sum $A=\displaystyle{\sum_{j \ge 0}}c_jz^jp_j \in  H^2\otimes L$ and for every $k \ge 0$, we have 
$$S^{*k}A=\displaystyle{\sum_{0\le j <k}}c_jS^{*(k-j)}p_j+\displaystyle{\sum_{ j\ge k}}c_jz^{j-k}p_j \in F+S^*L+\ldots+S^{*N}L.$$
(Note that $S^{*k}L=\{0\}$ for $k>N$). Then,
$$E_{H^2\otimes L} \subset clos(F+S^*L+\ldots+S^{*N}L) \subset \displaystyle{\bigcap_{k \ge 0}}E_{S^{*k}f}.$$
(For the last inclusion: $\displaystyle{\bigcap_{k \ge 0}}E_{S^{*k}f}$ is a closed $S^*$-invariant subspace and $F,\, L \subset \displaystyle{\bigcap_{k \ge 0}}E_{S^{*k}f}$.\\
According to $(ii)$, $f-p \in F=H^2\otimes L$ where $p$ is a polynomial. Then there exists $d \ge 0$ such that for every $k \ge d$, we have 
$$S^{*k}f \in S^{*d}(H^2\otimes L) \subset clos(F+S^*L+\ldots+S^{*N}L).$$ 
This leads to
$$\displaystyle{\bigcap_{k \ge 0}}E_{S^{*k}f} \subset clos(F+S^*L+\ldots+S^{*N}L),$$
and then holds the equality in $(iii)$.\\
$(iv)$ Observe first that in the case where $dim\, X<\infty$, the subspace $F$ has the finite co-dimension in $E_f$ (see the formula in $(iii)$ of the Theorem).
Suppose now that $E_f=H^2(\mathbb{D},\, X)$. Then, $H^2\otimes L$ is a $S$-invariant subspace in $H^2(\mathbb{D},\, X)$ of finite co-dimension. It is easy to see that such a subspace has the maximal local rank (for example, we can use the representation of Lax-Halmos which gives $F=H^2\otimes L=B H^2(\mathbb{D},\, X^\prime)$ (where $B$ is a product of Blaschke-Potapov). But, since $F= H^2\otimes L$, the local rank of $F$, i.e. 
$$r:=\displaystyle{\max_{|z|<1}}\, dim(f(z):\, f\in F),$$
is equal to the local rank of $L$. Then, $r=dim\,X$.\\
Conversely, suppose that $r=dim\,X$. The subspace $F\subset H^2(\mathbb{D},\, X)$ is $S$-invariant and then has the canonical representation of Lax-Halmos, $F=\Theta H^2(\mathbb{D},\, X^\prime)$, where $X^\prime \subset X$ and $\Theta$ is a left inner function. Since the local rank of $F$ coincides with the local rank of $L$ and using the Hypothesis, $r=dim\, X$, we have 
$$dim\, X=r=dim\, (\Theta(z) H^2(\mathbb{D},\, X^\prime)) \le dim\, X^\prime \le dim\, X.$$
Then, $X=X^\prime$ and,
$$det(\Theta)H^2(\mathbb{D},\, X )\subset \Theta H^2(\mathbb{D},\, X)\subset E_f.$$
But $\theta=det(\Theta)$ is a scalar inner function and the $S^*$-invariant subspace generated by $\theta H^2(\mathbb{D},\, X)$ coincides with $H^2(\mathbb{D},\, X)$. If $\phi_n$ are the Fej\'er polynomials of $\theta$, we have $$\displaystyle{\lim_{n \rightarrow \infty}}\|P_+\overline{\phi}_n\theta f-f\|_2=\displaystyle{\lim_{n \rightarrow \infty}}\|\phi_n(S^*)\theta f-f\|_2=0,$$ 
for every $f \in H^2(\mathbb{D},\, X)$ and therefore, $E_f=H^2(\mathbb{D},\, X)$.
 \hfill{$\Box$} 
\end{preuve}
\begin{rem}
If $dim\, X<\infty$, the subspace $H^2\otimes L$ has the finite co-dimension in $E_f$ and it seems very possible to have a certain criterion of cyclicity in $H^2(\mathbb{D},\, X)$ of lacunary series with bounded blocks.
\end{rem}
We deduce the following Corollary,
\begin{coro}
Let $f$ as in Theorem \ref{teobloc}.
\begin{enumerate}[(i)]
\item If $X=\mathbb{C}$, then $f$ is cyclic.
\item If $dim\, X<\infty$, then $H^2 \otimes L$ is a subspace of $E_f$ with the finite co-dimension ($\le (N+1).\, dim\, X$).
\item If $dim\, L=(N+1).\, dim\, X$, then $L=\mathcal{P}_N(X)$, and then $f$ is cyclic.
\item We have, $L\subset \displaystyle{\bigcap_{n \ge 1}} span_X(P_k:\, k \ge n)$ and 
$$\displaystyle{\bigcap_{n \ge 1}} span_X(P_k(z):\, |z|<1,\, k \ge n)=\displaystyle{\bigcap_{n \ge 1}} span_X(\widehat{P}_k(j):\, 0 \le j \le N,\, k \ge n).$$
\item If $\displaystyle{\bigcap_{n \ge 1}} span_X(P_k:\, k \ge n)\ne \mathcal{P}_N(X)$, or  $\displaystyle{\bigcap_{n \ge 1}} span_X(\widehat{P}_k(j):\, 0 \le j \le N,\, k \ge n)\ne X$, then $f$ is not cyclic in $H^2(\mathbb{D},\, X)$.
\end{enumerate}
\end{coro}
\textbf{Example:}
Let $X=\mathbb{C}^2,\, P_k=a_k(e_1+ze_2),\, a_k \ne 0,\, \displaystyle{\sum_{k \ge 0}}|a_k|^2<\infty$ o\`{u} $e_j,\, j=1,\, 2$ is the standard base of $\mathbb{C}^2$. Then, $N=1,\, L=\{\lambda P_1:\, \lambda \in \mathbb{C}\}$, the local rank of $L$ is $1$ and $\displaystyle{\bigcap_{n \ge 1}} span_X(P_k:\, k \ge n)\ne \mathcal{P}_1(\mathbb{C}^2)$. the function $f$, of course, is not cyclic (see Theorem \ref{zed}).
\section{The case of the polydisc}
In this section, we work with $H^2$ spaces on the polydisc $\mathbb{D}^n$ and consider multiparameter backward shifts. By definition,  
$$H^2(\mathbb{D}^n)=\{ f(z)=\displaystyle{\sum_{\alpha \ge 0}}\widehat{f}(\alpha)z^\alpha,\,\displaystyle{\sum_{\alpha \ge 0}}|\widehat{f}(\alpha)|^2<\infty\},$$
where $z=(z_1,\ldots,\, z_n) \in \mathbb{D}^n,\, \alpha=(\alpha_1,\ldots,\,\alpha_n) \in \mathbb{Z}_+^n$ is a multi-index and $z^\alpha=z_1^{\alpha_1}\ldots z_n^{\alpha_n}$ the elementary monomial in $\mathbb{D}^n$. We refer to $\cite{Rud1}$ for any supplementary information on the spaces $H^2(\mathbb{D}^n)$. We have $L^2(\mathbb{Z}_+^n)=l^2(\mathbb{Z}_+^n)\equiv H^2(\mathbb{D}^n)$.\\ 
We define the multiparameter semi-group $(S^{*\alpha})_{\alpha \in \mathbb{Z}_+^n}$ such that for any power series  $f(z)=\displaystyle{\sum_{\alpha \ge 0}}\widehat{f}(\alpha )z^\alpha$ in $\mathbb{D}^n$, we have
$$S^{*\alpha}f(z)=\displaystyle{\sum_{\beta \ge 0}}\widehat{f}(\alpha+\beta)z^\beta.$$
The semi-group $(S^{*\alpha})_{\alpha \in \mathbb{Z}_+^n}$ has $n$ generators $S_1^*,\ldots,\,S_n^*$. \\
It is clear that $f\in H^2(\mathbb{D}^n)\Rightarrow S^{*\alpha}f \in H^2(\mathbb{D}^n)$ and $\|S^{*\alpha}f\|_2 \le \|f\|_2$. Let $\sigma(f)$ be the (Fourier) spectrum of $f$,
$$\sigma(f)=\{ \alpha \in \mathbb{Z}_+^n: \, \widehat{f}(\alpha) \ne 0 \}.$$
We define the space $H^2$ on $\mathbb{D}^n$ with values in a separable Hilbert space $X$ by 
$$H^2(\mathbb{D}^n,\, X):=\{f(z)=\displaystyle{\sum_{\alpha \ge 0}}\widehat{f}(\alpha)z^\alpha \mbox{ such that } \widehat{f}(\alpha) \in X \mbox{ and } \displaystyle{\sum_{\alpha \ge 0}}\|\widehat{f}(\alpha)\|_X^2<\infty\}.$$
As before, we also define,
$$E_f=span_{H^2(\mathbb{D}^n,\, X)}\{S^{*\alpha}:\, \alpha \in \mathbb{Z}_+^n\}.$$
And $f$ is said cyclic if $E_f= H^2(\mathbb{D}^n,\, X)$.\\
The goal of this part is to find an analogue of the main result of the first part of this paper. We study the case of functions $f \in H^2(\mathbb{D}^n,\, X)$ whose spectrum is ``rare'' in the following sense.
\begin{enumerate}[(C1)]
\item There exists a constant $C$ such that,
$$card\{(\alpha,\, \alpha^\prime) \subset \sigma(f)\times \sigma(f),\, \alpha \ne \alpha^\prime:\, \beta=\alpha-\alpha^\prime \} \le C,\, \, \, \, \forall \, \, \beta \in \mathbb{Z}_+^n$$.
\item $\displaystyle{\lim_{j\rightarrow \infty}}(\alpha_{j+1,\, k}-\alpha_{j,\, k})=\infty,\, \, \, \, \forall\, \, k,\, 1 \le k \le n$.
\end{enumerate}
\textbf{Remark:}
It is obvious that if for a sequence $(\alpha_j)_{j\ge 1} \subset  \mathbb{Z}_+^n$, one of his componant sequence $(\alpha_{j,\, k})_{j \in \mathbb{Z}_+^n},\, 1 \le k \le n$ is a lacunary sequence in the sense of Hadamard, then the statement $(C1)$ is satisfied.
\begin{lem}\label{4poly}
Let $X$ a separable Hilbert space.\\
For any series $f \in H^2(\mathbb{D}^n,\, X)$ satisfying $(C1),\, (C2)$ above and if $(\widehat{f}(k))_{k \ge 0}$ is a relatively compact sequence, there exists a non-zero element $x \in X$ such that $$H^2(\mathbb{D}^n,\, X)\otimes x \subset E_{{S}^{* \beta}f} \, \, \, \, \forall  \, \beta \in \mathbb{Z}_+^n.$$
\end{lem}
\begin{preuve}
Let $\beta=(\beta_1,\ldots,\beta_n) \in \mathbb{Z}_+^n$ be fixed and $f\in H^2(\mathbb{D}^n,\, X)$ verifying the previous statements. According to $(C2)$, $f$ can be written  
$$f(z)= \displaystyle{\sum_{k \ge 0}} \widehat{f}(\alpha_k)z^{\alpha_k},$$
Using $(C2)$, there exists an integer $k_0$ such that $$\alpha_k -\alpha_{k-1} \ge \beta,\, \, \, \, \forall \, \, k \ge k_0.$$
And for any $k \ge k_0$,
$$\frac{S^{\alpha_k -\beta}}{\| \alpha_k \|} f(z)=\frac{\widehat{f}(\alpha_k)}{\| \widehat{f}(\alpha_k) \|}z^\beta+z^\beta \displaystyle{\sum_{j >k}} \frac{\widehat{f}(\alpha_j)}{\| \widehat{f}(\alpha_k )\|} \, z^{\alpha_j-\alpha_k}.$$
We take $r_k=\displaystyle \sum_{j >k} \frac{\widehat{f}(\alpha_j)}{\| \widehat{f}(\alpha_k) \|} \,
z_1^{\alpha_j-\alpha_k}$ and using the same method as in Lemma \ref{4} with some minor changements we prove that, 
$$ z^\beta x \in E_f \, \, \, \forall \,\, \beta \in \mathbb{Z}_+^n.$$
By the same way, we show that $ z^\beta x \in E_{S^{*\alpha}f} \, \, \, \forall \,\, \alpha,\, \beta \in \mathbb{Z}_+^n$ and $H^2(\mathbb{D}^n,\, X)\otimes x \subset  E_{S^{\alpha}f}\, \, \, \forall \,\, \alpha\in \mathbb{Z}_+^n .$
\hfill{$\Box$}
\end{preuve}
\begin{theo}\label{poly}
Let $X$ be a separable Hilbert space.\\
And $f\in H^2(\mathbb{D}^n,\, X),\, f(z)= \displaystyle{\sum_{k \ge 0}} \widehat{f}(\alpha_k) z^{\alpha_k}$ such that $\sigma(f)$ is satisfying statements $(C1),\, (C2)$\\ and $(\widehat{f}(\alpha_k))_{k \ge 0}$ is c.r.c. The following assertions are equivalents.
\begin{enumerate}[1)]
\item $f$ is cyclic for $S^*$ in $H^2(\mathbb{D}^n,\, X).$
\item $span(\widehat{f}({\alpha_j}) \, : \, j \ge k)=X \, \, \, \forall \, \, k \ge 0$.
\end{enumerate}
\end{theo}
\begin{preuve}
The study we give in the first part of this article with the space $H^2(\mathbb{D},\, X)$ where $X$ is a separable Hilbert space and the techniques developped can be used for the space $H^2(\mathbb{D}^n,\, X)$ since Lemma \ref{4poly} is fullfiled.
\hfill{$\Box$}
\end{preuve}
\begin{lem} $\cite{Abak}$\label{lemaba}
Let $1 \le p < \infty$, and let $\Omega$ be a $S^*$-invariant subspace of $\ell_a^p$ (i.e., if $f \in \Omega$ then $S^*f \in \Omega$). Suppose that the inclusion $1 \in span(S^{*k}f, \, k \ge 0)$ holds for any $f \in \Omega$. Then all elements of $\Omega$ are cyclic vectors in $\ell_a^p$.
\end{lem}
In the case $X\equiv \mathbb{C}$, we can easily prove the cyclicity of $f$, it is possible with Lemma \ref{4poly} to have $1 \in span(S^{*\alpha}f, \, \alpha \in \mathbb{Z}_+^2) $ because it is possible to divide by $x$ in the scalar case. To conclude, we need to give a generalization of the previous Lemma for the space $H^2(\mathbb{D}^n,\, X)$. If we consider $S^{*\alpha}=S_1^{*\alpha_1}\ldots S_n^{*\alpha_2} \, \, \forall \, \, \alpha=(\alpha_1,\ldots,\, \alpha_2) \in \mathbb{Z}_+^n$, it suffices to use the same induction argument done by E. Abakumov in Lemma \ref{lemaba} but in $H^2(\mathbb{D}^n)$ and in the following way; $\eta=(\eta_1,\ldots,\, \eta_n) \in \mathbb{Z}_+^n$ is fixed and we suppose that $z^\beta \in span(S^{*\alpha}f, \, \alpha \in \mathbb{Z}_+^n) \, \, \, \, \forall \, \, \beta < \eta$ then $z^\eta \in span(S^{*\alpha}f, \, \alpha \in \mathbb{Z}_+^n).$\\
\\
\begin{rem} It is interesting to see that in Theorem \ref{poly}, it is possible to take series $f(z)= \displaystyle{\sum_{k \ge 0}} \widehat{f}(\alpha_k) z^{\alpha_k}$ in the space $H^2(\mathbb{D}^n,\, X)$ whose Fourier spectrum is a set of the form $$\alpha_k=(m_k^1,\ldots,\, m_k^n),$$ 
where $(m_k^j)_{k \ge 1}$ are any lacunary sequences for $j=1,\ldots n$.
\hfill{$\blacksquare$}
\end{rem}

\end{document}